\documentclass[12pt]{amsart}
 
\usepackage[dvips, bookmarks=false, hyperfigures=false,
setpagesize=false]{hyperref}

\usepackage{graphicx}

\usepackage{pifont}

\usepackage[mathscr]{eucal}
\usepackage{t-angles}

\DeclareMathAlphabet{\mathpzc}{OT1}{pzc}{m}{it}

\advance\textheight\topskip
\usepackage{mathabx}
\usepackage{mathptmx}

\usepackage[matrix,arrow,curve,dvips]{xy} 
\UsePSspecials{dvips}

\allowdisplaybreaks

\numberwithin{equation}{section}
\numberwithin{figure}{section}
\makeatletter
\@addtoreset{equation}{section}
\@addtoreset{figure}{section}
\@addtoreset{subsubsection}{section}

\def\@secnumfont{\bfseries}
\def\subsubsection{\@startsection{subsubsection}{3}%
  \z@{.5\linespacing\@plus.7\linespacing}{-.5em}%
  {\normalfont\bfseries}}
\def\paragraph{\@startsection{paragraph}{4}%
  \z@\z@{-\fontdimen2\font}%
  \normalfont\bfseries}
\def\subparagraph{\@startsection{subparagraph}{5}%
  \z@\z@{-\fontdimen2\font}%
  \normalfont\bfseries}
\makeatother

\usepackage[OT2,T1]{fontenc}
\newcommand{\mifody}{%
  \renewcommand\rmdefault{wncyr}%
  \renewcommand\sfdefault{wncyss}%
  \renewcommand\encodingdefault{OT2}%
  \normalfont
  \selectfont}
\usepackage{times}

\newcommand{\Ribbon}{\boldsymbol{\pmb{\vartheta}}}

\newcommand{\eigen}[1]{\lambda(#1)}

\newcommand{\q}{\mathfrak{q}}

\newcommand{\coeff}[6]{c^{#1\kern-1pt,\,#2}_{#3\kern-1pt,\,#4}(#5\kern-1pt,#6)}

\newcommand{\VertexI}[2]{V^{#2}_{#1}}
\newcommand{\UertexI}[2]{U^{#2}_{#1}}

\newcommand{\VertexIII}[4]{V^{#2,\,#4}_{#1,\,#3}}
\newcommand{\UertexIII}[4]{U^{#2,\,#4}_{#1,\,#3}}

\newcommand{\VertexV}[6]{V^{#2,\,#4,\,#6}_{#1,\,#3,\,#5}}

\newcommand{\VertexM}[4]{V^{#1,\,#2}_{#3,\,#4}}
\newcommand{\UertexM}[4]{U^{#1,\,#2}_{#3,\,#4}}

\newcommand{\actsright}{\mathbin{\mbox{\large$\looparrowdownleft$}}}

\newcommand{\dualX}{{}^{^{\vee}}\kern-2pt\modX}
\newcommand{\dualXi}[1]{{}^{^{\vee}}\kern-2pt\modX^{\{#1\}}_0}

\newcommand{\RepXi}[1]{\modX^{\{#1\}}_0}
\newcommand{\RepXii}[4]{\modX^{\{#2,\,#4\}}_{#1,\,#3}}
\newcommand{\RepSii}[4]{\modS^{\{#2,\,#4\}}_{#1,\,#3}}
\newcommand{\RepPii}[4]{\modP^{\{#2,\,#4\}}_{#1,\,#3}}
\newcommand{\RepVii}[4]{\modV^{\{#2,\,#4\}}_{#1,\,#3}}
\newcommand{\RepLii}[4]{\modL^{\{#2,\,#4\}}_{#1,\,#3}}
\newcommand{\RepBii}[4]{\modB^{\{#2,\,#4\}}_{#1,\,#3}}

\newcommand{\RepA}{\modA}
\newcommand{\RepB}{\modB}
\newcommand{\RepL}{\modL}
\newcommand{\RepS}{\modS}
\newcommand{\RepX}{\modX}
\newcommand{\RepP}{\modP}
\newcommand{\RepV}{\modV}
\newcommand{\RepW}{\modW}
\newcommand{\RepY}{\modY}
\newcommand{\RepZ}{\modZ}

\newcommand{\DualW}{{}^{\vee}\!\RepW}
\newcommand{\DualZ}{{}^{\vee}\!\RepZ}

\newcommand{\fX}{\,\boldsymbol{\mathsf{x}}}
\newcommand{\fP}{\,\boldsymbol{\mathsf{p}}}

\newcommand{\Mod}[2]{#1\;\text{mod}\;#2}

\newcommand{\modP}{\mathscr{P}}
\newcommand{\modV}{\mathscr{V}}
\newcommand{\modW}{\mathscr{W}}
\newcommand{\modY}{\mathscr{Y}}
\newcommand{\modZ}{\mathscr{Z}}

\newcommand{\modA}{\mathscr{A}}
\newcommand{\modX}{\mathscr{X}}
\newcommand{\modL}{\mathscr{L}}
\newcommand{\modB}{\mathscr{B}}
\newcommand{\modS}{\mathscr{S}}

\newcommand{\Pdot}{\Pi_{_{\bullet}}}

\newcommand{\OPi}{\O{\Pdot}}

\newcommand{\Times}{\mbox{\LARGE${\times}$}}
\newcommand{\TIMES}{\mbox{\LARGE${\bigtimes}$}}

\newcommand{\Nich}{\mathfrak{B}}

\newcommand{\leftact}{\kern1pt{\rightharpoonup}\kern1pt}
\newcommand{\rightact}{\kern1pt{\leftharpoonup}\kern1pt}
\newcommand{\leftreg}{\kern1pt{\rightharpoonup}\kern1pt}
\newcommand{\rightreg}{\kern1pt{\leftharpoonup}\kern1pt}
\newcommand{\rightregX}{\kern1pt{\leftbarharpoon}\kern1pt}
\newcommand{\leftregii}{\kern2pt{\looparrowright}\kern1pt}
\newcommand{\rightregii}{\kern1pt{\looparrowleft}\kern1pt}

\newcommand{\leftregU}{\kern1pt{\rightharpoondown}\kern1pt}
\newcommand{\rightregU}{\kern1pt{\leftharpoondown}\kern1pt}
\newcommand{\rightregUX}{\kern1pt{\barleftharpoon}\kern1pt}
\newcommand{\leftregUii}{\kern2pt{\looparrowdownright}\kern1pt}
\newcommand{\rightregUii}{\kern1pt{\looparrowdownleft}\kern2pt}

\newcommand{\eval}[2]{\bigl\langle#1,\,#2\bigr\rangle\,}

\newcommand{\oZ}{\mathbb{Z}}

\newcommand{\abin}[2]{\mathchoice%
  {\abinomm{#1}{#2}}{\abinommm{#1}{#2}}%
  {\abinommm{#1}{#2}}{\abinommm{#1}{#2}}}
\newcommand{\abinomm}[2]{\mbox{\footnotesize$\displaystyle
    \genfrac{\langle}{\rangle}{0pt}{}{#1}{#2}$}}
\newcommand{\abinommm}[2]{\genfrac{\langle}{\rangle}{0pt}{}{#1}{#2}}

\newcommand{\Aint}[1]{\langle#1\rangle}
\newcommand{\Afac}[1]{\langle#1\rangle!\,}
\newcommand{\Abin}[2]{\mathchoice%
  {\Abinm{#1}{#2}}{\Abinmm{#1}{#2}}%
  {\Abinmm{#1}{#2}}{\Abinmm{#1}{#2}}}
\newcommand{\Abinm}[2]{\mbox{\footnotesize$\displaystyle
    \genfrac{\langle}{\rangle}{0pt}{}{#1}{#2}$}}
\newcommand{\Abinmm}[2]{\genfrac{\langle}{\rangle}{0pt}{}{#1}{#2}}

\newcommand{\thev}[1]{F(#1)}

\newcommand{\alphaminus}{-\sqrt{\kern-1pt\frac{2}{p}}}

\newcommand{\KK}{\mathcal{K}\kern-5.7pt\raisebox{-3.9pt}{\footnotesize\textit{2}}\,}

\newcommand{\cross}{\textstyle\!\!{\times}\!\!}
\newcommand{\punct}{\textstyle{\circ}}

\newcommand{\dotact}{\mathbin{\pmb{.}}}

\newcommand{\id}{\mathrm{id}}
\newcommand{\tensor}{\otimes}

\newcommand{\ccirc}{\mathbin{\raisebox{1pt}{\,$\scriptscriptstyle\circ$\,}}}

\newcommand{\A}{\raisebox{.5pt}{\large$\mathpzc{S}\kern-1pt$}}
\newcommand{\medA}{\mathpzc{S}}
\newcommand{\hA}{\mbox{\large$\mathpzc{s}\kern-.8pt$}}

\newcommand{\YDname}{\mathcal{Y\kern-3ptD}}
\newcommand{\BByd}{{}\mbox{\small${}^{\Nich}_{\Nich}$}\YDname}

\newcommand{\Shuffle}{\mathop{\text{\mifody\sf Sh}}\nolimits}

\newcommand{\Shift}[2]{#2^{\uparrow#1}}

\newcommand{\shift}{\uparrow}

\newcommand{\Bbin}[2]{{\Shuffle^{#1}_{#2}}}
\newcommand{\Bfac}[1]{\mathfrak{S}_{#1}}

\newcommand{\fobject}[1]{\object{\scriptstyle #1}}

\newcommand{\adjoint}{%
  \mathchoice{\mathbin{\blacktriangleright}}%
  {\mathbin{\mbox{\small${\blacktriangleright}$}}}%
  {\mathbin{{\blacktriangleright}}}%
  {\mathbin{{\blacktriangleright}}}}

\newcommand{\Bbraid}{\pmb{\Psi}}

\newcommand{\ffrac}[2]{\raisebox{.5pt}{\mbox{\footnotesize$\displaystyle\frac{#1}{#2}$}}}

\newcommand{\half}{%
  \mathchoice{\ffrac{1}{2}}{\frac{1}{2}}{\frac{1}{2}}{\frac{1}{2}}}

\newcommand{\cY}{\mathscr{Y}}

\newcommand{\bref}[1]{\textbf{\ref{#1}}}

\swapnumbers

\theoremstyle{definition}

\begin{document}

\title[Fusion for a rank-1 Nichols algebra]{Fusion in the entwined
  category of Yetter--Drinfeld modules of a rank-1 Nichols algebra}

\author[Semikhatov]{A.M.~Semikhatov}

\address{Lebedev Physics Institute, Russian Academy of Sciences}

\begin{abstract}
  We rederive a popular nonsemisimple fusion algebra in the braided
  context, from a Nichols algebra.  Together with the decomposition
  that we find for the product of simple Yetter--Drinfeld modules,
  this strongly suggests that the relevant Nichols algebra furnishes
  an equivalence with the triplet $W$-algebra in the $(p,1)$
  logarithmic models of conformal field theory.  For this, the
  category of Yetter--Drinfeld modules is to be regarded as an
  \textit{entwined} category (the one with monodromy, but not with
  braiding).
\end{abstract}

\maketitle

\section{Introduction}\thispagestyle{empty}%
The idea to construct ``purely algebraic'' counterparts of
vertex-operator algebras (conformal field theories) has a relatively
long history~\cite{[KLx],[Fink],[T],[MS],[FRS],[FFRS]}.
In~\cite{[FGST],[FGST2], [FGST3],[FGSTq],[NT],[GT],[BFGT],[BGT]}, this
idea was developed for nonsemisimple---logarithmic---CFT{} models,
which have been intensively studied recently
(see~\cite{[AA],[CR],[HY],[AN],[VJS],[AM-2p],[HLZ],[GRW],
  [FSS-11],[FSS-12],[VGJS],[GV],[RGW-12]} and the references therein).
In~\cite{[STbr]}, further, a braided and arguably ``more fundamental''
algebraic counterpart of logarithmic CFT{} was proposed.  It is given
by Nichols
algebras~\cite{[Nich],[AG],[AS-onthe],[AS-pointed],[Andr-remarks]};
the impressive recent progress in their theory (see
\cite{[Heck-class],[Heck-Weyl],[AHS],[ARS],[GHV],[GH-lyndon],[AFGV],
  [AAY],[Ag-0804-standard],[Ag-1008-presentation],[Ag-1104-diagonal]}
and the references therein) is a remarkable ``spin-off'' of
Andruskiewitsch and Schneider's program of classification of pointed
Hopf algebras.

Associating Nichols algebras with CFT{} models implies that certain
CFT-related structures must be reproducible from (some) Nichols
algebras.  Here, we take the simplest, rank-1 Nichols algebra
$\Nich\!_p$ of dimension $p\geq2$ and, from the category of its
Yetter--Drinfeld modules, extract a commutative associative
$2p$-dimensional algebra on the $\fX(r)_{\nu}$, $1\leq r\leq p$,
$\nu\in\oZ_2$:
\begin{equation}\label{the-algebra}
  \fX(r_1)_{\nu_1}\fX(r_2)_{\nu_2}
  = \sum_{\substack{s=|r_1-r_2|+1\\
      \text{step}=2}}^{p-1-|r_1+r_2-p|}
  \fX(s)_{\nu_1+\nu_2}
  \,\,+\,\,\sum_{\substack{s = 2p-r_1-r_2+1\\ \text{step}=2}}^{p}
  \fP(s)_{\nu_1+\nu_2},
\end{equation}
with
\begin{equation*}
  \fP(r)_{\nu}=
  \begin{cases}
    2\fX(r)_{\nu} + 2\fX(p-r)_{\nu+1},& r<p,\\
    \fX(p)_\nu,& r=p.
  \end{cases}
\end{equation*}
This is the FHST fusion algebra~\cite{[FHST]} (also see~\cite{[GT]}),
which makes part of what we know from \cite{[NT]} (also
see~\cite{[AM-3]}) to be an equivalence of representation
categories---of the \textit{triplet algebra} $W(p)$ in the $(p,1)$
logarithmic conformal
models~\cite{[Kausch],[Gaberdiel-K],[Gaberdiel-K-2],
  [Gaberdiel-K-3],[FHST]} and of a small quantum $s\ell_2$ at the
$2p$th root of unity, proposed in this capacity
in~\cite{[FGST],[FGST2]} and then used and studied, in particular,
in~\cite{[MN],[FHT],[Ar],[KoSa],[S-yd]} (this quantum group had
appeared before in~\cite{[AGL],[Su],[X]}).

The reoccurrence of the fusion algebra in the braided approach
advocated in~\cite{[STbr]}, together with some other observations,
supports the idea that Nichols algebras are \textit{at least as good
  as} the quantum groups proposed
previously~\cite{[FGST],[FGST2],[FGST3],[FGSTq],[S-q]} for the
logarithmic version of the Kazhdan--Lusztig correspondence (the
correspondence between vertex-operator algebras and quantum
groups).\footnote{But the actual motivation in~\cite{[STbr]}, which is
  yet to be tested on more advanced examples, was that Nichols
  algebras can actually do better than the ``old'' quantum groups.}

Algebra~\eqref{the-algebra} arises here as an algebra in the center of
the category of Yetter--Drinfeld $\Nich\!_p$-modules; the $\fX(r)_{\nu}$
are certain images of the simple Yetter--Drinfeld $\Nich\!_p$-modules
$\RepX(r)_{\nu}$.\footnote{The notation is fully explained below, but
  here we note that the module comodule structure, e.g., of
  $\RepX(r)_\nu$ depends only on $r$, whereas $\nu$ serves to
  distinguish isomorphic module comodules that nevertheless have
  different braiding.}  More is actually true: from the study of the
representation theory of~$\Nich\!_p$, we obtain that the tensor product
of simple Yetter--Drinfeld $\Nich\!_p$-modules decomposes as
\begin{equation}\label{the-fusion}
  \RepX(r_1)_{\nu_1}\tensor\RepX(r_2)_{\nu_2}
  = \bigoplus_{\substack{s=|r_1-r_2|+1\\
      \text{step}=2}}^{p-1-|r_1+r_2-p|}
  \RepX(s)_{\nu_1+\nu_2}
  \,\,\oplus\,\,\bigoplus_{\substack{s = 2p-r_1-r_2+1\\ \text{step}=2}}^{p}
  \RepP[s]_{\nu_1+\nu_2},
\end{equation}
where $\RepP[p]_{\nu}=\RepX(p)_{\nu}$ and $\RepP[r]_{\nu}$ for $1\leq
r\leq p-1$ is a reducible Yetter--Drinfeld $\Nich\!_p$-module with the
structure of subquotients
\begin{equation}\label{P-mod-ind}
  \raisebox{-2\baselineskip}{$\RepP[r]_{\nu} \ =\
    \ $}
  \xymatrix@=12pt{
    &\modX(p-r)_{\nu+1}
    \ar@/_12pt/[dl]
    \ar[dr]
    &\\
    \modX(r)_{\nu}\ar[dr] &&\modX(r)_{\nu+2},\ar@/_12pt/[dl]
    \\
    &\modX(p-r)_{\nu+1}&
  }
\end{equation}
Decompositions~\eqref{the-fusion} were conjectured in~\cite{[STbr]}
and are proved here.  The $\RepX(r)_\nu$ and $\RepP[r]_\nu$ do not
exhaust all the category of Yetter--Drinfeld $\Nich\!_p$-modules, but
make up ``the most significant part of it,'' and
relations~\eqref{the-fusion}, together with the structure of
$\RepP[r]_\nu$, already seem to imply that the category of
Yetter--Drinfeld $\Nich\!_p$-modules is equivalent to the $W(p)$
representation category.  This requires an important clarification,
however.

In the braided category of Yetter--Drinfeld $\Nich\!_p$-modules, the
simple objects are the $\RepX(r)_{\nu}$ labeled by $1\leq r\leq p$ and
$\nu\in\oZ_4$ (and, accordingly, $\nu\in\oZ_4$ in $\RepP[r]_\nu$, and
so on).  There are twice as many objects as in the category of $W(p)$
representations~\cite{[FHST],[AM-3],[NT]}.  But the presumed
equivalence is maintained for \textit{entwined}
categories~\cite{[Brug]}---those endowed with only ``double braiding''
$D_{\RepY,\RepZ}=c_{\RepZ,\RepY}\ccirc c_{\RepY,\RepZ}$ (the
\textit{monodromy} on the $W(p)$ side).  The properties of double
braiding can be axiomatized\pagebreak[3] without having to resort to
the braiding itself~\cite{[Brug]}.  This defines a \textit{twine
  structure} and, accordingly, an entwined category.  Remarkably, it
was noted in~\cite{[Brug]} that
\begin{quotation}
  ``many significant notions apparently related to $c$ actually depend
  only on $D$ or [the twist] $\theta$.  The $S$-matrix, and the
  subcategory of transparent objects, which play an important role in
  the construction of invariants of 3-manifolds, are defined purely in
  terms of the double braiding. More surprisingly, the invariants of
  ribbon links \dots\ do not depend on the actual braiding, but only
  on $D$.''
\end{quotation}
In the entwined category of Yetter--Drinfeld $\Nich\!_p$-modules, the
objects with $\nu$ and $\nu+2$ in their labels are isomorphic, which
sets $\nu\in\oZ_2$ and resolves the ``representation doubling
problem''; everything else on the algebraic side appears to be already
``fine-tuned'' to ensure the equivalence.  (We do not go as far as
modular transformations in this paper, but the above quotation
suggests that dealing with entwined categories is not an impediment to
rederiving the $W(p)$ modular properties at the Nichols algebra level,
in a ``braided version'' of what was done in~\cite{[FGST]}.)

It may also be worth noting that we derive~\eqref{the-algebra}
and~\eqref{the-fusion} independently (of course, from the same
structural results on Yetter--Drinfeld $\Nich\!_p$-modules, but not from
one another).  In particular, \eqref{the-algebra} is obtained by
directly composing the \textit{action} of $\fX(r_1)_{\nu_1}$ and
$\fX(r_2)_{\nu_2}$ on Yetter--Drinfeld modules, with
$\fX(r)_{\nu}:\RepY\to\RepY$ given by ``running $\RepX(r)_{\nu}$ along
the loop'' in the diagram (with the notation to be detailed in what
follows)
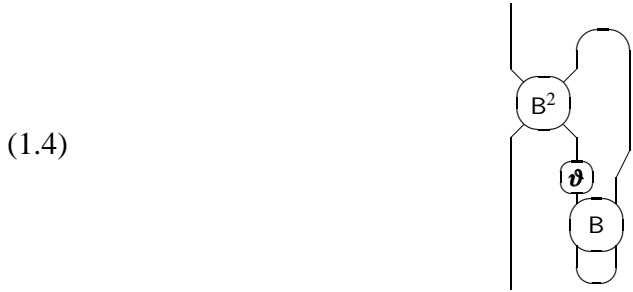
\begin{equation}\label{the-loop}
  \begin{tangles}{l}
    \id\step[2.5]\coev\\
    \vstr{50}\dh\step[1.5]\ddh\step[2]\id\\[-2pt]
    \step[.25]\obox{2}{\mathsf{B}^2}\step[2.25]\id\\[-2pt]
    \vstr{50}\hdd\step[1.5]\hd\step[2]\id\\
    \id\step[2.5]\O{\Ribbon}\step[1.5]\hdd\\[-2pt]
    \id\step[2.25]\obox{2}{\mathsf{B}}\\[-2pt]
    \vstr{33}\id\step[2.5]\id\step[1.5]\id\\
    \hh\id\step[2.75]{\makeatletter\@ev{0,\hm@de}{5,\hm@detens}{15}b\makeatother}
  \end{tangles}
\end{equation}

\noindent
As such, the $\fX(r)_{\nu}$ depend only on $\nu\in\oZ_2$---there is no
``$\oZ_4$ option'' for them.\footnote{Diagram~\eqref{the-loop}
  involves not only the squared braiding $\mathsf{B}^2$ of
  Yetter--Drinfeld modules but also, ``in the loop,'' the braiding
  itself (and the ribbon map~$\Ribbon$).  This does not affect the
  statement of the equivalence of entwined categories, but rather
  suggests exploring a further possibility, elaborating on the fact
  that the braiding of a Yetter--Drinfeld $\Nich\!_p$-module
  \textit{with itself} and \textit{with its dual} also depends on
  $\nu\in\oZ_2$, not $\nu\in\oZ_4$ (and the same for the ribbon map).
  An entwined$'$ category might allow these braidings in addition to
  twines. This is similar to the idea of \textit{twist equivalence} in
  the theory of Nichols algebras~\cite{[AS-pointed]} (the similarity
  is not necessarily superficial if we recall that the braiding of
  ``bare vertex operators'' is diagonal for $\Nich\!_p$).}

This paper is organized as follows.  For the convenience of the
reader, we summarize the relevant points from~\cite{[STbr]} in
Sec.~\ref{sec:Nich}; a very brief summary is that for a Nichols
algebra~$\Nich(X)$,\pagebreak[3] a category of its Yetter--Drinfeld
modules can be constructed using another braided vector space~$Y$
(whose elements are here called ``vertices,'' and the Yetter--Drinfeld
modules the ``multivertex'' modules).  In Sec.~\ref{sec:duality}, we
introduce duality and the related assumptions that make it possible to
write diagrams~\eqref{the-loop}.  In Sec.~\ref{sec:p}, everything is
specialized to a rank-1 Nichols algebra $\Nich\!_p$ (depending on an
integer $p\geq2$).  First and foremost, ``everything'' includes
multivertex Yetter--Drinfeld modules.  We actually construct important
classes of these modules quite explicitly
(Appendix~\bref{app:modules}), which allows proving~\eqref{the-fusion}
and also establishing duality relations among the modules.  We also
study their braiding, find the ribbon structure, and finally use all
this to derive~\eqref{the-algebra} from~\eqref{the-loop}
for~$\Nich\!_p$.  Basic properties of Yetter--Drinfeld modules over a
braided Hopf algebra are recalled in Appendix~\bref{app:YD-axiom}.

\section{The Nichols algebra of screenings}\label{sec:Nich}
We summarize the relevant points of~\cite{[STbr]} in this section.

\subsection*{Screenings and $\Nich(X)$}
The underlying idea is that the non\-locali\-ties associated with
screening operators---multiple-integration contours, such as
\begin{equation}\label{crosses}
  \xymatrix@R=4pt@C=80pt{
    \ar@{--}|(.2){\cross}|(.55){\cross}|(.8){\cross}[0,2]&&
  }
  = \iiint_{-\infty<z_1<z_2<z_3<\infty} 
  s_{i_1}(z_1) s_{i_2}(z_2) s_{i_3}(z_3),
\end{equation}
where $s_j(z)$ are the ``screening currents''--- allow introducing a
\textit{coproduct} by contour cutting, called ``deconcatenation'' in
what follows:
\begin{align}\label{cutting}
  \Delta:\xymatrix@R=4pt@C=26pt{
    \ar@{--}|(.25){\cross}|(.5){\cross}|(.75){\cross}[0,2]&&
  }\mapsto{}& \xymatrix@R=4pt@C=9pt{
    \ar@{--}|(.25){\cross}|(.5){\cross}|(.75){\cross}[0,3]
    &&&{\ \raisebox{-6pt}{\rotatebox{90}{\mbox{\ding{34}}}}\ }\ar@{--}[0,3]&&& } + \xymatrix@R=4pt@C=9pt{
    \ar@{--}|(.33){\cross}|(.66){\cross}[0,3]
    &&&{\ \raisebox{-6pt}{\rotatebox{90}{\mbox{\ding{34}}}}\ }\ar@{--}|(.5){\cross}[0,3]&&& }
  \\
  \notag &+ \xymatrix@R=4pt@C=9pt{ \ar@{--}|(.5){\cross}[0,3]
    &&&{\ \raisebox{-6pt}{\rotatebox{90}{\mbox{\ding{34}}}}\ }\ar@{--}|(.33){\cross}|(.66){\cross}[0,3]&&& } +
  \xymatrix@R=4pt@C=9pt{ \ar@{--}[0,3]
    &&&{\ \raisebox{-6pt}{\rotatebox{90}{\mbox{\ding{34}}}}\ }\ar@{--}|(.25){\cross}|(.5){\cross}|(.75){\cross}[0,3] &&& }
\end{align}
(with the line cutting symbol subsequently understood as~$\tensor$).
A \textit{product} of ``lines populated with crosses'' is also
defined, as the ``quantum'' shuffle product~\cite{[Rosso-CR]}, which
involves a \textit{braiding} between any two screenings.  It is well
known that these three structures---coproduct, product, and
braiding---satisfy the braided bialgebra axioms~\cite{[Rosso-CR]}.
The \textit{antipode} is in addition given by contour reversal.  The
braided Hopf algebra axioms are then satisfied for quite a general
braiding (by far more general than may be needed in CFT); it is rather
amusing to see how the braided Hopf algebra axioms are satisfied by
merging and cutting contour integrals~\cite{[STbr]}.  The algebra
\textit{generated by} single crosses---individual screenings---is the
Nichols algebra $\Nich(X)$ of the braided vector space $X$ spanned by
the different screening species (whose number is called the rank of
the Nichols algebra).  

\subsection*{Nichols algebras}
The Nichols algebras---``bialgebras of type one'' in
\cite{[Nich]}---are a crucial element in a classification program of
\textit{ordinary} Hopf algebras of a certain type (see
\cite{[AG],[AS-pointed],[AS-onthe],[ARS]} and the references therein).
Nichols algebras have several definitions, whose equivalence is due
to~\cite{[Sch-borel]} and~\cite{[AG]}.  The Nichols algebra $\Nich(X)$
of a braided linear space $X$ can be characterized as a graded braided
Hopf algebra $\Nich(X)=\bigoplus_{n\geq0}\Nich(X)^{(n)}$ such that
$\Nich(X)^{(1)}=X$ and this last space coincides with the space of
\textit{all primitive elements} $P(X)=\{x\in\Nich(X) \mid\Delta
x=x\tensor 1 + 1\tensor x\}$ and it \textit{generates all of
  $\Nich(X)$} as an algebra.\footnote{An important technicality, noted
  in \cite{[A-about],[G-free]}, is a distinction between quantum
  symmetric algebras~\cite{[Rosso-inv]} and Nichols algebras proper;
  the latter are selected by the condition that the braiding be
  \textit{rigid}, which in particular guarantees that the duals $X^*$
  are objects in the same braided category with the $X$.}  Nichols
algebras occurred independently in~\cite{[Wor]}, in constructing a
quantum differential calculus, as ``fully braided generalizations'' of
symmetric algebras,
\begin{equation*}
  \Nich(X)
  = k\oplus X\oplus\bigoplus_{r\ge 2} X^{\otimes r}/\ker\Bfac{r},
\end{equation*}
where $\Bfac{r}$ is the total braided symmetrizer (``braided
factorial'').

\subsection*{The space of vertices $Y$}
In addition to the braided linear space $X$ spanned by the different
screening species, we introduce the space of vertex operators taken at
a fixed point,
\begin{equation}
  Y=\text{Span}(V_\alpha(0)),
\end{equation}
where $\alpha$ ranges over the different primary fields in a given
CFT{} model.  CFT{} also yields the braiding $\Psi:X\tensor X\to
X\tensor X$ of any two screenings (which is always applied to two
screenings on the same line, as in~\eqref{crosses}), as well as the
braiding $\Psi:X\tensor Y\to Y\tensor X$ and $\Psi:Y\tensor X\to
X\tensor Y$ of a screening and vertex (also on the same line, as
in~\eqref{1-punctured} below), and eventually the braiding
$\Psi:Y\tensor Y\to Y\tensor Y$ of any two vertices, but a large part
of our construction can be formulated without this last.

The two braided vector spaces $X$ and $Y$ are all that we need in this
section; the braiding $\Psi$ can be entirely general.

\subsection*{Dressed vertex operators as $\Nich(X)$-modules}
We use the space $Y$ to construct $\Nich(X)$-modules.  Their elements
are sometimes referred to in CFT{} as ``dressed$/$screened vertex
operators,'' for example,
\begin{equation}\label{1-punctured}
  \xymatrix@R=4pt@C=70pt{
    \ar@*{[|(1.6)]}@{-}|(.1){\cross}|(.35){\cross}|(.55){\punct}|(.7){\cross}[0,2]&&\\
  }
  =
  \iint_{-\infty<x_1<x_2<0}\!\!\!\!\!\!\!
  s_{i_1}(x_1) s_{i_2}(x_2)\; V_\alpha(0)\!\!\!
  \int\limits_{0<x_3<\infty} s_{i_3}(x_3).
\end{equation}
It is understood that the $\times$ and $\circ$ are decorated with the
appropriate indices read off from the right-hand side; but it is in
fact quite useful to suppress the indices altogether and let $\times$
and $\circ$ respectively denote the entire spaces $X$ and $Y$, and we
assume this in what follows.

Because the integrations can be taken both on the left and on the
right of the vertex position, the resulting modules are actually
$\Nich(X)$ bimodules.  The left and right actions of $\Nich(X)$ are by
pushing the ``new'' crosses into the different positions using
braiding; the left action, for example, can be visualized as
\begin{multline*}
  \xymatrix@R=6pt@C=30pt{
    \ar@{--}|(.5){\cross}[0,2]&&
  }
  \dotact
  \xymatrix@R=6pt@C=30pt{
    \ar@*{[|(1.6)]}@{-}|(.3){\punct}|(.65){\cross}[0,2]&&
  }
  =
  \\[6pt]
  \xymatrix@R=6pt@C=30pt{
    \ar@*{[|(1.6)]}@{-}|(.1){\cross}|(.3){\punct}|(.70){\cross}[0,2]&&
  } + 
  \xymatrix@R=6pt@C=30pt{
    \ar@*{[|(1.6)]}@{-}|(.3){\punct}|(.5){\cross}|(.75){\cross}[0,2]
    \ar@/^10pt/[r]&&
  } + 
  \xymatrix@R=6pt@C=30pt{
    \ar@*{[|(1.6)]}@{-}|(.3){\punct}|(.6){\cross}|(.95){\cross}[0,2]
    \ar@/^12pt/[rr]&&
  }
\end{multline*}
where the arrows, somewhat conventionally, represent the
braiding~$\Psi$.  Once again by deconcatenation, e.g.,
\begin{align*}
  \delta_{\text{L}}:\xymatrix@R=4pt@C=28pt{
    \ar@*{[|(1.6)]}@{-}|(.25){\cross}|(.5){\cross}|(.65){\punct}|
    (.85){\cross}[0,2]&&}
  \mapsto{}&
  \xymatrix@R=4pt@C=12pt{ \ar@{--}[0,3]
    &&&{\ \raisebox{-6pt}{\rotatebox{90}{\mbox{\ding{34}}}}\ }\ar@*{[|(1.6)]}@{-}|(.25){\cross}|(.5){\cross}|(.65){\punct}|(.85){\cross}[0,3] &&&}
  +
  \xymatrix@R=4pt@C=12pt{ \ar@{--}|(.5){\cross}[0,3]
    &&&{\ \raisebox{-6pt}{\rotatebox{90}{\mbox{\ding{34}}}}\ }\ar@*{[|(1.6)]}@{-}|(.3){\cross}|(.65){\punct}|(.85){\cross}[0,3]
    &&&}
  \\
  &{}+\xymatrix@R=4pt@C=12pt{\ar@{--}|(.3){\cross}|(.6){\cross}[0,3]
    &&&{\ \raisebox{-6pt}{\rotatebox{90}{\mbox{\ding{34}}}}\ }\ar@*{[|(1.6)]}@{-}|(.65){\punct}|(.85){\cross}[0,3]
    &&&},
\end{align*}
these bimodules are also bicomodules and, in fact, Hopf bimodules
over~$\Nich(X)$ (see
\cite{[Besp-TMF],[Besp-next],[Besp-Dr-(Bi)],[Majid-book]} for the
general definitions).

\subsection*{Braid group diagrams and quantum shuffles}
A standard graphical representation for the multiplication in
$\Nich(X)$ and its action on its modules is in terms of braid group
diagrams.  For example, the above left action is represented as (to be
read from top down)
\begin{equation}\label{Sh12}
  \begin{tangles}{l}
    \fobject{{\times}}\step\fobject{\otimes}\step
    \fobject{{\circ}}\step\fobject{{\times}}\\
    \vstr{200}\id\step[2]\id\step\id
  \end{tangles}\
  \to
  \begin{tangles}{l}
    \fobject{{\times}}\step
    \fobject{{\circ}}\step\fobject{{\times}}\\
    \vstr{200}\id\step[1]\id\step\id
  \end{tangles}\
  + \
  \begin{tangles}{l}
    \fobject{{\times}}\step
    \fobject{{\circ}}\step[1]\fobject{{\times}}\\
    \vstr{200}\hx\step\id
  \end{tangles}\
  + \ \ 
  \begin{tangles}{l}
    \fobject{{\times}}\step
    \fobject{{\circ}}\step[1]\fobject{{\times}}\\
    \hx\step\id\\
    \id\step\hx
  \end{tangles}
  \ \ = (\id + \Psi_{1} + \Psi_{2}\Psi_{1})(X\tensor Y\tensor X),
\end{equation}

\smallskip

\noindent
where we use the ``leg notation,'' in the right-hand side, letting
$\Psi_i$ denote the braiding of the $i$th and $(i+1)$th factors in a
tensor product (our notation and conventions are the same as
in~\cite{[STbr]}).  The braid group algebra element
$\Bbin{}{1,2}\equiv\id + \Psi_{1} + \Psi_{2}\Psi_{1}$ occurring here
is an example of quantum shuffles.  The product in $\Nich(X)$ is in
fact the shuffle product
\begin{equation}\label{Sh-prod}
  \Bbin{}{r,s}:X^{\otimes r}\tensor X^{\otimes s}\to X^{\otimes(r+s)}
\end{equation}
on each graded subspace.  The antipode restricted to each $X^{\otimes
  r}$ is up to a sign given by the ``half-twist''---\,the braid group
element obtained via the Matsumoto section from the longest element in
the symmetric group:
\begin{equation}\label{antipode}
  \A_r=
  (-1)^r\,
  \Psi_1 (\Psi_2\Psi_1)(\Psi_3\Psi_2\Psi_1)\dots
  (\Psi_{r-1}\Psi_{r-2}\dots\Psi_1):X^{\otimes r}
  \to X^{\otimes r}
\end{equation}
(with the brackets inserted to highlight the structure, and the sign
inherited from reversing the integrations); for
example,\enlargethispage{\baselineskip}
\begin{equation*}
  \A_5=-
  \ \begin{tangles}{l}
 \hstr{70}\vstr{50}\hx\step[1]\id\step[1]\id\step[1]\id\\
 \hstr{70}\vstr{50}\id\step[1]\hx\step[1]\id\step[1]\id\\
 \hstr{70}\vstr{50}\id\step[1]\id\step[1]\hx\step[1]\id\\
 \hstr{70}\vstr{50}\hx\step[1]\id\step[1]\hx\\
 \hstr{70}\vstr{50}\id\step[1]\hx\step[1]\id\step[1]\id\\
 \hstr{70}\vstr{50}\hx\step[1]\hx\step[1]\id\\
 \hstr{70}\vstr{50}\id\step[1]\hx\step[1]\id\step[1]\id\\
 \hstr{70}\vstr{50}\hx\step[1]\id\step[1]\id\step[1]\id\\
 \end{tangles}\ 
\end{equation*}
The Hopf bimodules alluded to above are (some subspaces in)
$\bigoplus\limits_{r,s\geq 0}X^{\otimes r}\tensor Y\tensor X^{\otimes
  s}$, with the left and right $\Nich(X)$ actions on these also
expressed in terms of quantum shuffles as
\begin{align*}
  \Bbin{}{r,s+1+t}&:X^{\otimes r}\tensor
  \bigl(X^{\otimes s}\tensor Y\tensor X^{\otimes t}\bigr)\to
  \bigoplus_{i=0}^r X^{\otimes(s+r-i)}\tensor Y\tensor
  X^{\otimes(t+i)}\pagebreak[3]\\[-6pt]
  \intertext{and}
  \Bbin{}{s+1+t,r}&:
  \bigl(X^{\otimes s}\tensor Y\tensor X^{\otimes t}\bigr)
  \tensor X^{\otimes r}\to
  \bigoplus_{i=0}^r X^{\otimes(s+r-i)}\tensor Y\tensor
  X^{\otimes(t+i)}.
\end{align*}

\subsection*{Hopf-algebra diagrams}
The four operations on bi(co)modules of a braided Hopf algebra
$\Nich$ are standardly expressed as
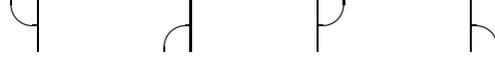
\begin{equation*}
  \begin{tangles}{ccccccc}
    \lu&\qquad\qquad&\ld&\qquad\qquad&\ru&\qquad\qquad&\rd\\
  \end{tangles}
\end{equation*}

\smallskip

\noindent
which are respectively the left module structure
$\Nich\tensor\RepZ\to\RepZ$, the left comodule
$\RepZ\to\Nich\tensor\RepZ$, the right module structure
$\RepZ\tensor\Nich\to\RepZ$, and the right comodule structures
$\RepZ\to\RepZ\tensor\Nich$.  The product and coproduct in the braided
Hopf algebra itself are denoted as
$\begin{tangles}{l}
  \hcu
\end{tangles}$ \ and \
$\begin{tangles}{l}
  \hcd
\end{tangles}$\;.  The braiding is still denoted as
$\begin{tangles}{l} \vstr{67}\hx
\end{tangles}\ $, but in contrast to the braid-group diagrams, each
line now represents a copy of $\Nich$ or a $\Nich$ (co)module.

\subsection*{Adjoint action and Yetter--Drinfeld modules}
The left and right actions of a braided Hopf algebra $\Nich$ on its
Hopf bimodule $\RepZ$ give rise to the \textit{left adjoint action}
$\Nich\tensor\RepZ\to\RepZ$:
\begin{equation}\label{adja}
  \begin{tangles}{l}
    \vstr{240}\lu\object{\raisebox{22pt}{\kern-4pt\tiny$\blacktriangleright$}}
  \end{tangles}
  \ \ = \ \
  \begin{tangles}{l}
    \hcd\step\id\\
    \vstr{80}\id\step\hx\\
    \lu\step\O{\medA}\\
    \vstr{75}\step\ru
  \end{tangles}
\end{equation}
A fundamental fact is that \textit{the space of right coinvariants in
  a Hopf bimodule is invariant under the left adjoint action}; this
actually leads to an equivalence of categories, the category of Hopf
bimodules and the category of Yetter--Drinfeld
modules~\cite{[Besp-TMF],[Besp-next],[Sch-H-YD],[Wor]}.  We recall
some relevant facts about Yetter--Drinfeld modules in
Appendix~\ref{app:YD-axiom}.  In our case of modules spanned by
dressed vertex operators, the right coinvariants---all those $y$ that
map as $y\mapsto y\tensor 1$ under the right coaction---are simply the
vertex operators dressed by screenings only from the left, i.e.,
elements of $X^{\otimes r}\tensor Y$, for example,
$\xymatrix@C=40pt@1{
  \ar@*{[|(1.6)]}@{-}|(.3){\cross}|(.6){\cross}|(.75){\punct}[0,2]&&
}$.
In terms of \textit{braid group} diagrams (with the lines representing
the $X$ and $Y$ spaces), an example of the left adjoint action on such
spaces is given by
\begin{equation}\label{adja.1.2}
  \begin{tangles}{l}
    \hstr{80}\fobject{{\times}}\step\fobject{\otimes}\step
    \fobject{{\times}}\step\fobject{{\times}}\step\fobject{{\circ}}\\
    \hstr{80}\vstr{280}\id\step[2]\id\step\id\step\id
  \end{tangles}\
  \to\
  \begin{tangles}{l}
    \\
    \hstr{80}\vstr{320}\id\step\id\step\id\step\id
  \end{tangles}
  \;+
  \
  \begin{tangles}{l}
    \hstr{80}\vstr{50}\id\step\id\step\id\step\id\\
    \hstr{80}\vstr{220}\hx\step\id\step\id\\
    \hstr{80}\vstr{50}\id\step\id\step\id\step\id\\
  \end{tangles}
  \;+
  \;
  \begin{tangles}{l}
    \hstr{80}\vstr{160}\hx\step\id\step\id\\
    \hstr{80}\vstr{160}\id\step\hx\step\id
  \end{tangles}
  \;-
  \;
  \begin{tangles}{l}\vstr{90}
    \hstr{90}\vstr{80}\hx\step\id\step\id\\
    \hstr{90}\vstr{80}\id\step\hx\step\id\\
    \hstr{90}\vstr{80}\id\step\id\step\hx\\
    \hstr{90}\vstr{80}\id\step\id\step\hx
  \end{tangles}
  \;-
  \;
  \begin{tangles}{l}
    \hstr{90}\vstr{66}\hx\step\id\step\id\\
    \hstr{90}\vstr{66}\id\step\hx\step\id\\
    \hstr{90}\vstr{66}\id\step\id\step\hx\\
    \hstr{90}\vstr{66}\id\step\id\step\hx\\
    \hstr{90}\vstr{66}\id\step\hx\step\id
  \end{tangles}
  \;-
  \;
  \begin{tangles}{l}
    \hstr{90}\vstr{55}\hx\step\id\step\id\\
    \hstr{90}\vstr{55}\id\step\hx\step\id\\
    \hstr{90}\vstr{55}\id\step\id\step\hx\\
    \hstr{90}\vstr{55}\id\step\id\step\hx\\
    \hstr{90}\vstr{55}\id\step\hx\step\id\\
    \hstr{90}\vstr{55}\hx\step\id\step\id
  \end{tangles}
\end{equation}

\smallskip

\noindent
where a single ``new'' cross arrives to each of the three possible
positions in two ways, one with the plus and the other with the minus
sign in front (which is something expected of an ``adjoint'' action).
That the cross never stays to the right of $\circ$ is precisely a
manifestation of the above invariance statement for the space of right
coinvariants.  This means that a number of terms\pagebreak[3] that
follow when expressing~\eqref{adja} in terms of braid group diagrams
cancel.  The left adjoint action~\eqref{adja} can in fact be expressed
more economically as follows.

We define a modified left action $\ \begin{tangles}{l}\vstr{70}\lu
  \object{\raisebox{5.2pt}{\tiny$\bullet$}}
\end{tangles}\ $ of $\Nich(X)$ on its Hopf bimodules spanned by
dressed vertex operators by allowing the ``new'' crosses to arrive
only to the left of $\circ$, for example,
\begin{equation}\label{left-dot-example}
  \begin{tangles}{l}
    \hstr{80}\fobject{{\times}}\step\fobject{\otimes}\step
    \fobject{{\times}}\step\fobject{{\times}}\step\fobject{{\circ}}\\
    \hstr{80}\vstr{200}\id\step[2]\id\step\id\step\id
  \end{tangles}\
  \to\
  \begin{tangles}{l}
    \\
    \hstr{80}\vstr{240}\id\step\id\step\id\step\id
  \end{tangles}
  \;+
  \
  \begin{tangles}{l}
    \hstr{80}\vstr{240}\hx\step\id\step\id\\
  \end{tangles}
  \;+
  \;
  \begin{tangles}{l}
    \hstr{80}\vstr{120}\hx\step\id\step\id\\
    \hstr{80}\vstr{120}\id\step\hx\step\id
  \end{tangles}
\end{equation}

\noindent
(more crosses might be initially placed to the right of the vertex
$\circ$; the action does not see them).  In general, 
$\ \begin{tangles}{l}\vstr{70}\lu
  \object{\raisebox{5.2pt}{\tiny$\bullet$}}
\end{tangles}\ $ is the map
\begin{equation}\label{FromLeftii}
  \ \begin{tangles}{l}\vstr{70}\lu
  \object{\raisebox{5.2pt}{\tiny$\bullet$}}
\end{tangles}\  = \Bbin{}{r,s}:X^{\otimes r}\tensor
  \bigl(X^{\otimes s}\tensor Y\bigr)
  \to X^{\otimes(r+s)}\tensor Y.
\end{equation}
Similarly, a modified right action $\ \begin{tangles}{l}
  \vstr{70}\object{\raisebox{5.1pt}{\tiny$\bullet$}}\ru
\end{tangles}\ $ on the space of right coinvariants is defined by
first letting the new cross to be braided with the vertex and then
shuffling into all possible positions relative to the ``old'' crosses:
\begin{equation*}
  \begin{tangles}{l}
    \hstr{80}\fobject{{\times}}\step\fobject{{\times}}\step
    \fobject{{\circ}}\step\fobject{\otimes}\step\fobject{{\times}}\\
    \hstr{80}\vstr{200}\id\step\id\step\id\step[2]\id
  \end{tangles}\
  \to\
  \begin{tangles}{l}
    \hstr{80}\vstr{240}\id\step\id\step\hx\\
  \end{tangles}\ \
  + \ \
  \begin{tangles}{l}
    \hstr{80}\vstr{120}\id\step\id\step\hx\\
    \hstr{80}\vstr{120}\id\step\hx\step\id\\
  \end{tangles}\ \
  + \ \
  \begin{tangles}{l}
    \hstr{80}\vstr{80}\id\step\id\step\hx\\
    \hstr{80}\vstr{80}\id\step\hx\step\id\\
    \hstr{80}\vstr{80}\hx\step\id\step\id
  \end{tangles}
\end{equation*}

\smallskip

\noindent
which in general is
\begin{equation}\label{FromRightii}
  \ \begin{tangles}{l}
    \vstr{70}\object{\raisebox{5.1pt}{\tiny$\bullet$}}\ru
  \end{tangles}\ = \Bbin{}{s,r}\ccirc(\id^{\otimes s}\tensor\Bbraid_{1,r})
  :\bigl(X^{\otimes s}\tensor Y\bigr)\tensor X^{\otimes r} \to
  X^{\otimes(s+r)}\tensor Y,
\end{equation}
where $\Bbraid_{s,r}$ is the braiding of an $s$-fold tensor product
with an $r$-fold tensor product.  The
$\ \begin{tangles}{l}\vstr{70}\lu
  \object{\raisebox{5.2pt}{\tiny$\bullet$}}
\end{tangles}\ $ and $\ \begin{tangles}{l}
  \vstr{70}\object{\raisebox{5.1pt}{\tiny$\bullet$}}\ru
\end{tangles}\ $ actions preserve the spaces of right coinvariants and
commute with each other.  The ``economic'' expression for adjoint
action~\eqref{adja} is~\cite{[STbr]}
\begin{equation}\label{dot-adja}
  \begin{tangles}{l}
    \vstr{200}\lu
    \object{\raisebox{18pt}{\kern-4pt\tiny$\blacktriangleright$}}
  \end{tangles}
  \ \ = \ \
  \begin{tangles}{l}
    \hcd\step\id\\
    \vstr{80}\id\step\hx\\
    \vstr{100}\lu\object{\raisebox{8.4pt}{\tiny$\bullet$}}
    \step\O{\medA}\\
    \vstr{75}\step\object{\raisebox{5.8pt}{\tiny$\bullet$}}\ru
  \end{tangles}
\end{equation}
This diagram is the map
\begin{equation}\label{adja-formula}
  {\adjoint}_{r,s}\equiv
  \sum_{i=0}^{r} \Bbin{}{r - i, s + i} \bigl(\Bbin{}{s, i}
  \Shift{s}{\Bbraid_{1,i}} \,\Shift{(s+1)}{\A_i}\, \Bbraid_{i,s+1}
  \bigr)^{\shift(r-i)}:
  X^{\otimes r}\tensor(X^{\otimes s}\tensor Y)\to X^{\otimes(s+r)}\tensor Y.
\end{equation}

\subsection*{Multivertex Yetter--Drinfeld modules}
More general, \textit{multivertex}, Yetter--Drinfeld
$\Nich(X)$-modules can be constructed by letting two or more vertices
(the $Y$ spaces) sit on the same line, e.g.,
\begin{equation}\label{ex1}
  \xymatrix@R=6pt@C=50pt{
    \ar@*{[|(1.6)]}@{-}|(.15){\cross}
    |(.3){\punct}|(.45){\cross}|(.6){\cross}|(.75){\cross}|(.90){\punct}[0,2]&&
  }\quad\text{or}\quad
  \xymatrix@R=6pt@C=50pt{
    \ar@*{[|(1.6)]}@{-}|(.15){\cross}
    |(.25){\punct}|(.4){\cross}|(.5){\cross}|(.65){\punct}
    |(.8){\cross}
    |(.90){\punct}[0,2]&&
  }
\end{equation}
These diagrams respectively represent $X\tensor Y\tensor X^{\otimes
  3}\tensor Y$ and $X\tensor Y\tensor X^{\otimes 2}\tensor Y\tensor
X\tensor Y$ 
(in general, different spaces could be taken instead of
copies of the same~$Y$, but in our setting they are all the same).  By
definition, the $\Nich(X)$ action and coaction on these are
\begin{enumerate}\addtolength{\itemsep}{4pt}
  \addtocounter{equation}{1}
  \renewcommand{\labelenumi}{(\theequation.\arabic{enumi})}
  \renewcommand{\theenumi}{\thesection.\arabic{equation}}
\item\label{2items} the ``cumulative'' left adjoint action, and

\item deconcatenation up to the first $\circ$.
\end{enumerate}
The ``cumulative'' adjoint means that all the $\circ$ except the
rightmost one are viewed on equal footing with the $\times$ under this
action: the adjoint action of $X^{\otimes r}$ on the space $X^{\otimes
  s}\tensor Y\tensor X^{\otimes t}\tensor Y$ in a two-vertex module is
given by ${\adjoint}_{r,s+1+t}$.  
For example, the left adjoint action
$\xymatrix@1@C=30pt{ \ar@{--} |(.5){\cross}[0,2]&& } \ \adjoint \
\xymatrix@1@C=40pt{ \ar@*{[|(1.6)]}@{-}
  |(.3){\punct}|(.65){\cross}|(.90){\punct}[0,2]&& }$ is given by the
braid group diagrams that are exactly those in the right-hand side
of~\eqref{adja.1.2}, with the corresponding strand representing not
${\times}{}=X$ but~${\circ}{}=Y$.
The $\Nich(X)$ coaction by deconcatenation up to the first vertex
means, for example, that at most one $\times$ can be deconcatenated in
each diagram in~\eqref{ex1}.

For multivertex Yetter--Drinfeld modules, the form~\eqref{dot-adja} of
the adjoint action is valid if 
$\ \begin{tangles}{l}\vstr{70}\lu
  \object{\raisebox{5.2pt}{\tiny$\bullet$}}
\end{tangles}\ $ is understood as the ``cumulative'' action preserving
right coinvariants; for example,
\begin{equation*}
  \xymatrix@1@C=30pt{
    \ar@{--} |(.5){\cross}[0,2]&&
  }
  \ \mbox{\small$\bullet$} \ 
  \xymatrix@1@C=40pt{
    \ar@*{[|(1.6)]}@{-} |(.3){\punct}|(.65){\cross}|(.90){\punct}[0,2]&&
  }
\end{equation*}
is given just by the braid group diagrams in the right-hand side
of~\eqref{left-dot-example} with the second strand representing not
${\times}{}=X$ but~${\circ}{}=Y$.

\subsection*{Fusion product}
The multivertex Yetter--Drinfeld modules are not exactly tensor
products of single-vertex ones---they carry a different action, which
is \textit{not} $(\mu_{\RepY}\tensor\mu_{\RepZ})\ccirc\Delta$, and the
coaction is not diagonal either.  They actually follow via a
\textit{fusion product}~\cite{[STbr]}, which is defined on two
single-vertex Yetter--Drinfeld modules (each of which is the space of
right coinvariants in a Hopf bimodule) as
\begin{equation}\label{fusion-def}
  \begin{tangles}{l}
    \id\step[1]\ld\\
    \object{\raisebox{8pt}{\tiny$\bullet$}}\ru\step[1]\id
  \end{tangles}
\end{equation}

\noindent
which is the map
\begin{equation*}
  \sum_{j=0}^{t}
  \Bbin{}{s, j} \Shift{s}{\Bbraid_{1, j}}:
  (X^{\otimes s}\tensor Y)\tensor
  (X^{\otimes t}\tensor Y)\to
  X^{\otimes s}\tensor Y\tensor X^{\otimes t}\tensor Y
\end{equation*}
on each $(s,t)$ component.  
For example, if $s=2$ and $t=3$, the top
of the above diagram can be represented as
\begin{equation*}
  \xymatrix@1@C=50pt{
    \ar@*{[|(1.6)]}@{-}|(.15){\cross}
    |(.55){\cross}|(.80){\punct}[0,2]&&
  }\ \ \tensor \ \ 
  \xymatrix@1@C=50pt{
    \ar@*{[|(1.6)]}@{-}
    |(.25){\cross}|(.5){\cross}|(.7){\cross}
    |(.85){\punct}[0,2]&&
  }
\end{equation*}

\noindent
and then in view of the definition of \ $\ \begin{tangles}{l}
  \vstr{70}\object{\raisebox{5.1pt}{\tiny$\bullet$}}\ru
\end{tangles}\ $, the meaning of~\eqref{fusion-def} is that $j\geq 0$
crosses from the right factor are detached from their ``native''
module and sent to mix with the left crosses (the sum over $j$ is
taken in accordance with the definition of the coaction).

The construction extends by taking the fusion product of multivertex
modules: the coaction in~\eqref{fusion-def} is then the one just
described, by deconcatenation up to the first vertex, and the
$\ \begin{tangles}{l}
  \vstr{70}\object{\raisebox{5.1pt}{\tiny$\bullet$}}\ru
\end{tangles}\ $ action on a multivertex module is ``cumulative,''
i.e., each cross acting from the right, e.g., on $\xymatrix@1@C=50pt{
  \ar@*{[|(1.6)]}@{-}|(.15){\cross}
  |(.3){\punct}|(.5){\cross}|(.7){\cross} |(.90){\punct}[0,2]&& }$,
arrives at each of the \textit{five} possible positions.

\section{Duality in the category of Yetter--Drinfeld
  modules}\label{sec:duality}
We now consider duality in a braided category of representations of a
braided Hopf algebra $\Nich$.  We briefly recall the standard
definitions and basic properties, and then assume that duality exists
in the setting of the preceding section; this then allows us to
construct endomorphisms of the identity functor in Sec.~\ref{sec:p}.

\subsection{}
For a $\Nich$-module $\RepZ$, we let $\DualZ$ denote the left dual
module in the same (rigid) braided category.  The duality means that
there are coevaluation and evaluation maps
\begin{equation*}
  \begin{tangles}{l}
    \vphantom{x}\\
     \Coev\\
    \fobject{\RepZ}\fobject{\kern60pt\DualZ}
  \end{tangles}
  \qquad\qquad\text{and}\qquad
  \begin{tangles}{l}
    \fobject{\DualZ}\fobject{\kern60pt \RepZ}\\
    \Ev
  \end{tangles}
\end{equation*}

\noindent
which are morphisms in the category and satisfy the axioms
\begin{equation*}
  \begin{tangles}{l}
    \vstr{15}\id\\
    \vstr{67}\hstr{50}\id\step[2]\coev\\
    \vstr{67}\hstr{50}\ev\step[2]\id\\
    \vstr{67}\hstr{50}\vstr{33}\step[4]\id
  \end{tangles}
  \quad=\ \ \
  \begin{tangles}{l}
    \vstr{175}\id
  \end{tangles}
  \qquad\text{and}\qquad
  \begin{tangles}{l}
    \vstr{67}\hstr{50}\vstr{33}\step[4]\id\\
    \vstr{67}\hstr{50}\coev\step[2]\id\\
    \vstr{67}\hstr{50}\id\step[2]\ev\\
    \vstr{15}\id
  \end{tangles}
  \quad=\ \ \
  \begin{tangles}{l}
    \vstr{175}\id
  \end{tangles}
\end{equation*}

\smallskip

\noindent
where the two straight lines are $\id_{\DualZ}$ and $\id_{\RepZ}$. \ It
follows that
\begin{equation*}
  \begin{tangles}{l}
    \hxx\step[1]\id\\
    \hh\id\step\ev    
  \end{tangles}
  \ \ = \ \
  \begin{tangles}{l}
    \id\step[1]\hx\\
    \hh\ev\step[1]\id
  \end{tangles}
\end{equation*}

\noindent
and similarly for the coevaluation.

The dual $\DualZ$ to a left--left Yetter--Drinfeld $\Nich$-module
$\RepZ$ is a left--left Yetter--Drinfeld $\Nich$-module with the
action and coaction, \textit{temporarily} denoted by
$\ \begin{tangles}{l}\vstr{80}\lu \object{\raisebox{5.2pt}{$\oleft$}}
\end{tangles}\ $ \ and  $\ \begin{tangles}{l}\vstr{80}
  \ld\object{\raisebox{5.2pt}{$\oleft$}}
\end{tangles}\ $\;, \ defined as~\cite{[Besp-next]}
\begin{equation}\label{on-dual}
  \begin{tangles}{l}
    \id\step[1]\id\\
    \vstr{120}\lu[1]
    \object{\raisebox{9.3pt}{$\oleft$}}\\
    \step[1]\id
  \end{tangles}\ \ = \ \ 
  \begin{tangles}{l}
    \vphantom{x}\\
    \O{\medA}\step\id\step\coev\\
    \vstr{50}\hx\step[1]\id\step[2]\id\\
    \id\step[1]
    \lu\object{\raisebox{8pt}{\kern-4pt\tiny$\blacktriangleright$}}\step[2]\id\\
    \ev\step[2]\id
  \end{tangles}
  \qquad\text{and}\qquad
  \begin{tangles}{l}
    \step[1]\id\\
    \vstr{120}\ld[1]
    \object{\raisebox{9.3pt}{$\oleft$}}\\
    \id\step[1]\id
  \end{tangles}\ \ = \ \ 
  \begin{tangles}{l}
    \id\step[3]\coev\\
    \id\step[2]\ld\step[2]\id\\
    \vstr{50}\id\step[2]\hxx\step[2]\id\\
    \ev\step\O{\medA^{{\scriptscriptstyle-1}}}\step[2]\id
  \end{tangles}
\end{equation}
The definitions are equivalent to the properties (which, inter alia,
imply that the evaluation is a $\Nich$ module comodule morphism)
\begin{equation}\label{equiv-to}
  \begin{tangles}{l}
    \id\step[1]\id\step[2]\id\\
    \vstr{120}\lu[1]
    \object{\raisebox{9.3pt}{$\oleft$}}\step[2]\id
    \\
    \step[1]\ev
  \end{tangles}\ \ = \ \
  \begin{tangles}{l}
    \O{\medA}\step\id\step[1]\id\\
    \vstr{50}\hx\step[1]\id\\
    \id\step\lu\object{\raisebox{8pt}{\kern-4pt\tiny$\blacktriangleright$}}\\
    \ev
  \end{tangles}
  \qquad\text{and}\qquad
  \begin{tangles}{l}
    \vstr{50}\step[1]\id\step[2]\id\\
    \vstr{120}\ld[1]\object{\raisebox{9.3pt}{$\oleft$}}\step[2]\id\\[-4pt]
    \id\step[1]\ev
  \end{tangles}\ \ = \ \ 
  \begin{tangles}{l}
    \id\step[2]\ld\\
    \vstr{50}\id\step[2]\hxx\\
    \ev\step\O{\medA^{{\scriptscriptstyle-1}}}
  \end{tangles}
\end{equation}

We prove the Yetter--Drinfeld property for for
$\ \begin{tangles}{l}\vstr{80}\lu \object{\raisebox{5.2pt}{$\oleft$}}
\end{tangles}\ $ \ and  $\ \begin{tangles}{l}\vstr{80}
  \ld\object{\raisebox{5.2pt}{$\oleft$}}
\end{tangles}\ $ \ for completeness.  In view of~\eqref{equiv-to}, it
is easiest to verify the Yetter--Drinfeld axiom by establishing that
\begin{equation}\label{to-prove}
  \begin{tangles}{l}
    \vstr{90}\cd\step\id\step[2]\id\\
    \vstr{50}\id\step[2]\hx\step[2]\id\\
    \vstr{90}\lu[2] \object{\raisebox{6.2pt}{$\oleft$}}
    \step\hd\step[1.5]\id\\
    \vstr{90}\ld[2]\object{\raisebox{6.2pt}{$\oleft$}}\step\ddh
    \step[1.5]\id\\
    \vstr{50}\id\step[2]\hx\step[2]\id\\
    \vstr{90}\cu\step\ev
  \end{tangles}\ \ = \ \
  \begin{tangles}{l}
    \cd\step\ld\object{\raisebox{6.9pt}{$\oleft$}}\step[2]\id\\[-4pt]
    \id\step[2]\hx\step\id\step[2]\id\\
    \cu\step\lu
    \object{\raisebox{7.2pt}{$\oleft$}}\step[2]\id\\
    \step\id\step[3]\ev
  \end{tangles}
\end{equation}
Pushing the new action and then the coaction ``to the other side,'' we
see that the left-hand side of~\eqref{to-prove}, by the above
properties, is equal to
\begin{equation*}
  \begin{tangles}{l}
    \hcd\fobject{\quad\checkmark}\step[2]\id\step[1]\id\\[-4pt]
    \O{\medA}\step[1]\x\step[1]\id\\
    \vstr{50}\hx\step[2]\hxx\\
    \id\step[1]\hd\step[1]\hld\step\id\\
    \dh\step[1]\hx\step[.5]\id\step\id\\[-.8pt]
    \step[.5]\hx\step
    \hlu\object{\raisebox{13.2pt}{\kern-4pt\tiny$\blacktriangleright$}}
    \step[.5]\ddh\\
    \step[.5]\O{\medA^{{\scriptscriptstyle-1}}}\step[.75]
    {\makeatletter\@ev{0,\hm@de}{10,\hm@detens}{15}b\makeatother}
    \step[1.25]\dd\fobject{\checkmark}\\
    \step[.5]\cu
  \end{tangles}
  \quad=\quad
  \begin{tangles}{l}
    \O{\medA}\step[1]\id\step[2.5]\id\\
    \hx\step[2.5]\id\\
    \id\step[.5]\hcd\step[1]\ld\\
    \id\step[.5]\id\step[1]\hx\step\id\\
    \id\step[.5]\hcu\step
    \lu\object{\raisebox{8.1pt}{\kern-4pt\tiny$\blacktriangleright$}}\\
    \hx\step[2]\ddh\\
    \O{\medA^{{\scriptscriptstyle-1}}}\step[1]\ev
  \end{tangles}
  \quad=\quad
  \begin{tangles}{l}
    \O{\medA}\step[1]\id\step[1.5]\id\\
    \vstr{50}\hx\step[1.5]\id\\
    \id\step[.5]\hcd\step\id\\
    \vstr{50}\id\step[.5]\hd\step[.5]\hx\\
    \id\step[1]\hlu[1]\object{\raisebox{12.9pt}{\kern-4pt\tiny$\blacktriangleright$}}\step[-.5]\hld
    \step[1]\id\\
    \vstr{50}\id\step[.5]\ddh\step[.5]\hx\\
    \id\step[.5]\hcu\step\id\\
    \vstr{50}\hx\step[1.5]\id\\
    \O{\medA^{{\scriptscriptstyle-1}}}
    \step[.75]{\makeatletter\@ev{0,\hm@de}{10,\hm@detens}{15}b\makeatother}
  \end{tangles}
  \quad=\quad
  \begin{tangles}{l}
    \hcd\step\id\step[1.5]\id\\
    \vstr{50}\id\step\hx\step[1.5]\id\\
    \id\step\id\step\O{\medA}\step[.5]\dd\\
    \id\step\id\step\hlu\object{\raisebox{13pt}{\kern-4pt\tiny$\blacktriangleright$}}\step[-.5]\hld\\
    \id\step\id\step\O{\medA^{{\scriptscriptstyle-1}}}\step[.5]\d\\
    \vstr{50}\id\step\hx\step[1.5]\id\\
    \hcu\step[1.75]{\makeatletter\@ev{0,\hm@de}{0,\hm@detens}{15}b\makeatother}
  \end{tangles}
\end{equation*}
In the first diagram, we insert $\A$ at the position of the upper
checkmark and $\A^{-1}$ into the same line, at the lower checkmark,
and use the properties of the antipode,
\begin{equation*}
 \begin{tangles}{l}
  \cd\\
  \O{\medA}\step[2]\O{\medA}
    \end{tangles}\ = \
    \begin{tangles}{l}
      \step\O{\medA}\\
      \cd\\
      \xx\\
    \end{tangles}
     \qquad \text{and} \quad
     \begin{tangles}{l}
      \O{\medA^{{\scriptscriptstyle-1}}}\step[2]\O{\medA^{{\scriptscriptstyle-1}}}\\
      \cu
    \end{tangles}\ = \
    \begin{tangles}{l}
      \x\\
      \cu\\
      \step\O{\medA^{{\scriptscriptstyle-1}}}
    \end{tangles} 
\end{equation*}

\noindent
This readily gives the second diagram above, where we further
recognize the right-hand side of the Yetter--Drinfeld axiom assumed
for the module.  After using it (the third diagram), and after another
application of the properties of $\A$ and $\A^{-1}$, we obtain the
fourth diagram, and it is immediate to see that it coincides with the
right-hand side of~\eqref{to-prove} also rewritten by pushing
$\ \begin{tangles}{l}\vstr{80}\lu \object{\raisebox{5.2pt}{$\oleft$}}
\end{tangles}\ $ \ and  $\ \begin{tangles}{l}\vstr{80}
  \ld\object{\raisebox{5.2pt}{$\oleft$}}
\end{tangles}\ $\;\ ``to the other side.''

\subsection{Assuming a rigid category} We further assume that the
category of $n$-vertex Yetter--Drinfeld $\Nich$-modules is rigid; this
means that the dual modules are modules in the same category---in our
case, multivertex Yetter--Drinfeld $\Nich$-modules, and the action and
coaction defined in~\eqref{on-dual} are just those
in~\eqref{2items}---and hence the evaluation map satisfies the
properties
\begin{equation}\label{equiv-to-2}
  \begin{tangles}{l}
    \id\step[1]\id\step[2]\id\\
    \vstr{120}\lu[1]
    \object{\raisebox{10pt}{\kern-4pt\tiny$\blacktriangleright$}}
    \step[2]\id\\
    \step[1]\ev
  \end{tangles}\ \ = \ \
  \begin{tangles}{l}
    \O{\medA}\step\id\step[1]\id\\
    \vstr{50}\hx\step[1]\id\\
    \id\step\lu\object{\raisebox{8pt}{\kern-4pt\tiny$\blacktriangleright$}}\\
    \ev
  \end{tangles}
  \qquad\text{and}\qquad
  \begin{tangles}{l}
    \vstr{50}\step[1]\id\step[2]\id\\
    \vstr{120}\ld[1]\step[2]\id\\
    \id\step[1]\ev
  \end{tangles}\ \ = \ \ 
  \begin{tangles}{l}
    \id\step[2]\ld\\
    \id\step[2]\hxx\\
    \ev\step\O{\medA^{{\scriptscriptstyle-1}}}
  \end{tangles}
\end{equation}
for any pair of Yetter--Drinfeld $\Nich$-modules.  Evidently, we then
also have
\begin{equation}\label{equiv-to-3}
  \begin{tangles}{l}
    \id\step\coev\\
    \vstr{67}\hx\step[1]\dd\\
    \id\step\lu[1]\object{\raisebox{8pt}{\kern-4pt\tiny$\blacktriangleright$}}
  \end{tangles}
  \ \ =\quad
  \begin{tangles}{l}
    \O{\medA}\step\coev\\
    \lu\object{\raisebox{8pt}{\kern-4pt\tiny$\blacktriangleright$}}\step[2]\id\\
    \vstr{50}\step\id\step[2]\id
  \end{tangles}
  \qquad\text{and}\qquad
  \begin{tangles}{l}
    \coev\\
    \id\step\ld\\
    \hx\step[1]\id
  \end{tangles}
  \ \ = \ \
  \begin{tangles}{l}
    \step[1]\coev\\
    \ld\step[2]\id\\
    \O{\medA^{{\scriptscriptstyle-1}}}\step\id\step[2]\id
  \end{tangles}
\end{equation}

\subsection{}If the category $\BByd$ of Yetter--Drinfeld
$\Nich$-modules is rigid, then for each $\RepZ\in\BByd$, there is a
morphism $\chi_{\RepZ}:\RepY\to\RepY$ for any $\RepY\in\BByd$, defined
as
\begin{equation}\label{YD-loop}
  \begin{tangles}{l}
    \id\fobject{\kern10pt\RepY}\step[3]\fobject{\kern-12pt\RepZ}\coev
    \fobject{\kern12pt\DualZ}\\[-4pt]
    \vstr{50}\dh\step[1.5]\dd\step[2]\id\\[-2pt]
    \step[.25]\obox{2}{\mathsf{B}^2}\step[2.75]\id\\[-2pt]
    \vstr{50}\hdd\step[1.5]\d\step[2]\id\\
    \id\step[3]\O{\vartheta}\step[1]\dd\\
    \id\step[2.5]\obox{2}{\mathsf{B}}\\
    \vstr{33}\id\step[3]\id\step\id\\
    \hh\id\step[3]\ev
  \end{tangles}
  \qquad = \qquad
  \begin{tangles}{l}
    \id\fobject{\kern10pt\RepY}\step[3]\fobject{\kern-12pt\RepZ}\coev
    \fobject{\kern12pt\DualZ}\\[-4pt]
    \vstr{50}\dh\step[1.5]\dd\step[2]\id\\[-2pt]
    \step[.25]\obox{2}{\mathsf{B}^2}\step[2.75]\id\\[-2pt]
    \vstr{50}\hdd\step[1.5]\d\step[2]\id\\
    \id\step[3]\O{\vartheta}\step[2]\id\\
    \id\step[3]\O{\Times}\step[1]\dd\\
    \vstr{50}\id\step[3]\hx\\    
    \id\step[3]{\makeatletter\@ev{0,\hm@de}{5,\hm@detens}{10}b\makeatother}
  \end{tangles}
\end{equation}
where $\mathsf{B}$ is defined in~\eqref{B+B2} and $\vartheta$ is any
$\Nich$ module comodule morphism.  In the second diagram, Bespalov's
``squared relative antipode''~\cite{[Besp-TMF]}
\begin{equation}\label{eq:sigma2}
  \sigma_2\equiv\ \ \
  \begin{tangles}{l}
    \vstr{50}\id\\
    \O{\Times}\\
    \vstr{50}\id
  \end{tangles}
  \ \ \ = \ \ \
  \begin{tangles}{l}
    \ld[2]\\
    \O{\medA^2}\step[2]\OPi\\
    \x\\
    \object{\raisebox{8pt}{\tiny$\bullet$}}\ru[2]
  \end{tangles}
  \ \ \ = \ \ \
\begin{tangles}{l}
    \ld\\
    \O{\medA}\step\id\\
    \lu\object{\raisebox{8pt}{\kern-4pt\tiny$\blacktriangleright$}}\\
  \end{tangles}
  \qquad\text{such that}\qquad
  \begin{tangles}{l}
    \vstr{150}\lu\object{\raisebox{12.9pt}{\kern-4pt\tiny$\blacktriangleright$}}\\[-4pt]
    \step[1]\O{\Times}\\
    \step[1]\id
  \end{tangles}
  \quad =\quad
  \begin{tangles}{l}
    \x\\
    \O{\Times}\step[2]\O{\medA^2}\\
    \x\\
    \lu[2]\object{\raisebox{8pt}{\kern-4pt\tiny$\blacktriangleright$}}
  \end{tangles}
\end{equation}

\noindent
(see~\cite{[Besp-next],[BKLT]} for its further properties and use)
occurs in view of~\eqref{equiv-to-2}.

That the map defined by~\eqref{YD-loop} is a $\Nich$ module comodule
morphism follows from the general argument that so are $\mathsf{B}$,
evaluation, and coevaluation (and $\theta$).  It is also instructive
to see this by diagram manipulation (temporarily writing\quad
$\begin{tangles}{l} \O{\theta}
\end{tangles}$\quad for
\quad $\begin{tangles}{l}
  \O{\vartheta}\\
  \O{\Times}
\end{tangles}$\quad for brevity):
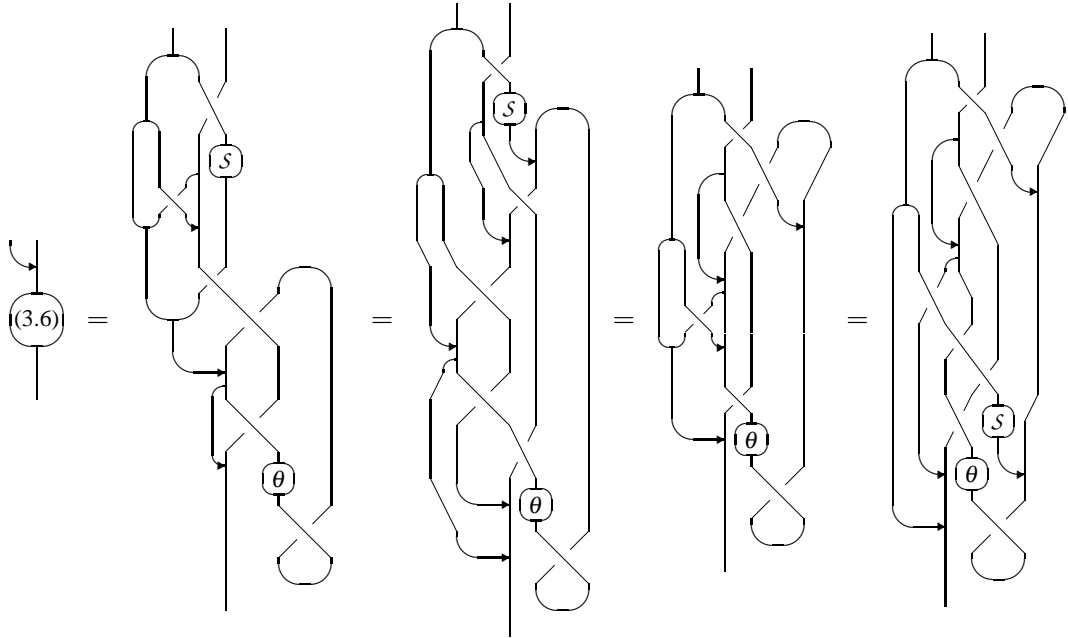
\begin{equation*}
  \begin{tangles}{l}
    \step[1]\lu\object{\raisebox{8pt}{\kern-4pt\tiny$\blacktriangleright$}}\\[-4pt]
    \step[1]\obox{2}{\eqref{YD-loop}}\\
    \step[2]\id
  \end{tangles}
  \ \ = \ \
  \begin{tangles}{l}
    \step[.5]\cd\step\id\\
    \vstr{100}\hcd\step[1.5]\hx\\
    \id\step[1]\id\step[1]\hld\step[1]\O{\medA}\\
    \vstr{50}\id\step[1]\hx\step[.5]\id\step\id\\[-.5pt]
    \hcu\step\hlu\object{\raisebox{13pt}{\kern-4pt\tiny$\blacktriangleright$}}\step\id\\[-10pt]
    \vstr{50}\step[.5]\id\step[2]\hx\step[2]\vstr{100}\coev\\
    \step[.5]\cu\step\x\step[2]\id\\
    \step[1.5]\lu[2]\object{\raisebox{8pt}{\kern-4pt\tiny$\blacktriangleright$}}\step[-.5]\hld\step[2]\id\step[2]\id\\
    \step[3]\id\step[.5]\x\step[2]\id\\[-.5pt]
    \step[3]\hlu\object{\raisebox{13pt}{\kern-4pt\tiny$\blacktriangleright$}}\step[2]
    \O{\theta}\step[2]\id\\
    \step[3.5]\id\step[2]\x\\
    \step[3.5]\id\step[2]\ev
  \end{tangles}
  \quad=\ \
  \begin{tangles}{l}
    \step[.5]\cd\step\id\\
    \vstr{50}\step[.5]\id\step[2]\hx\\
    \step[.5]\id\step[1.5]\hld\step[1]\O{\medA}\step[1]\coev\\
    \hcd\step[1]\dh\d\lu\object{\raisebox{8pt}{\kern-4pt\tiny$\blacktriangleright$}}\step[2]\id\\
    \vstr{50}\id\step[1]\id\step[1.5]\id\step[1]\hx\step[2]\id\\
    \dh\step[.5]\dh\step[1]\lu\object{\raisebox{8pt}{\kern-4pt\tiny$\blacktriangleright$}}\step[1]\id\step[2]\id\\
    \step[.5]\id\step[1]\x\step[1]\id\step[2]\id\\
    \step[.5]\lu\object{\raisebox{8pt}{\kern-4pt\tiny$\blacktriangleright$}}\step[-.5]\hld\step[2]\id\step[1]\id\step[2]\id\\
    \step[.5]\hdd\step[.5]\x\step[1]\id\step[2]\id\\
    \step[.5]\id\step[1]\id\step[2]\hx\step[2]\id\\
    \step[.5]\d\lu[2]\object{\raisebox{8pt}{\kern-4pt\tiny$\blacktriangleright$}}\step\O{\theta}\step[2]\id\\
    \step[1.5]\lu[2]\object{\raisebox{8pt}{\kern-4pt\tiny$\blacktriangleright$}}\step[1]\x\\
    \step[3.5]\id\step[1]\ev
  \end{tangles}
  \ \ =\ \
  \begin{tangles}{l}
    \step[.5]\cd\step\id\\[-10pt]
    \vstr{50}\step[.5]\id\step[2]\hx\step[1]\vstr{100}\coev\\
    \step[.5]\id\step[1]\ld\step[1]\hx\step[1]\dd\\
    \hcd\step[.5]\id\step[1]\hx\step[1]\lu[1]\object{\raisebox{8pt}{\kern-4pt\tiny$\blacktriangleright$}}\\
    \id\step[1]\id\step[.5]\lu\object{\raisebox{8pt}{\kern-4pt\tiny$\blacktriangleright$}}\step[-.5]\hld\step[1]\id\step[2]\id\\
    \vstr{50}\id\step[1]\hx\step[.5]\id\step\id\step[2]\id\\
    \hcu\step[1]\hlu\object{\raisebox{13pt}{\kern-4pt\tiny$\blacktriangleright$}}\step[1]\id\step[2]\id\\
    \vstr{50}\step[.5]\id\step[2]\hx\step[2]\id\\
    \step[.5]\lu[2]\object{\raisebox{8pt}{\kern-4pt\tiny$\blacktriangleright$}}\step[1]\O{\theta}\step[2]\id\\
    \step[2.5]\id\step[1]\x\\
    \step[2.5]\id\step[1]\ev
  \end{tangles}
  \ = \ \
  \begin{tangles}{l}
    \step[.5]\cd\step\id\\[-10pt]
    \vstr{50}\step[.5]\id\step[2]\hx\step[1]\vstr{100}\coev\\
    \step[.5]\id\step[1]\ld\step[1]\hx\step[1]\dd\\
    \hcd\step[.5]\id\step[1]\hx\step[1]\lu[1]\object{\raisebox{8pt}{\kern-4pt\tiny$\blacktriangleright$}}\\
    \id\step[1]\id\step[.5]\lu\object{\raisebox{8pt}{\kern-4pt\tiny$\blacktriangleright$}}\step[-.5]\hld\step[1]\hd\step[1.5]\id\\
    \id\step[1]\hx\step[.5]\hd\step[1]\id\step[1.5]\id\\
    \vstr{67}\id\step[1]\id\step[1]\hx\step[1]\id\step[1.5]\id\\
    \vstr{67}\id\step[1]\id\step[1]\id\step[1]\hx\step[1.5]\id\\
    \id\step[1]\id\step[1]\hx\step[1]\O{\medA}\step[1]\hdd\\
    \id\step\lu\object{\raisebox{8pt}{\kern-4pt\tiny$\blacktriangleright$}}\step[1]\O{\theta}\step[1]\lu[1]\object{\raisebox{8pt}{\kern-4pt\tiny$\blacktriangleright$}}\\
    \lu[2]\object{\raisebox{8pt}{\kern-4pt\tiny$\blacktriangleright$}}
    \step[1]\x\\
    \step[2]\id\step\ev
  \end{tangles}
\end{equation*}

\noindent
In the first equality, we use only the Yetter--Drinfeld axiom, with
$\mathsf{B}^2$ represented by the \textit{first} diagram for
$\mathsf{B}^2$ in~\eqref{B+B2}; the associativity of action was used
in the second equality above; another use of the associativity in the
lower part of the third diagram allows recognizing the left-hand side
of~\eqref{yd-axiom}; the Yetter--Drinfeld property is then applied in
the third equality together with the first property
in~\eqref{equiv-to-3}, yielding the fourth diagram; there, we use that
the property of $\sigma_2$ in~\eqref{eq:sigma2}
\noindent
and the first property in~\eqref{equiv-to-2} to obtain the last, fifth
diagram, where an ``antipode bubble'' is annihilated, showing that,
indeed,
\begin{equation*}
  \begin{tangles}{l}
    \step[1]\lu\object{\raisebox{8pt}{\kern-4pt\tiny$\blacktriangleright$}}\\[-4pt]
    \step[1]\obox{2}{\eqref{YD-loop}}\\
    \step[2]\id
  \end{tangles}
  \ \ = \ \
  \begin{tangles}{l}
    \id\step[2]\id\\
    \id\step[1]\obox{2}{\eqref{YD-loop}}\\
    \lu[2]\object{\raisebox{8pt}{\kern-4pt\tiny$\blacktriangleright$}}
  \end{tangles}
\end{equation*}

\noindent
The commutativity of~\eqref{YD-loop} with coaction can be verified
similarly.\enlargethispage{\baselineskip}

\subsection{Ribbon structure}A ribbon structure is a morphism
$\Ribbon:\RepY\to\RepY$ for every object $\RepY$ such that
\begin{equation}\label{Rib-def}
  \begin{tangles}{l}
  \step[.25]\O{\Ribbon}\step[1.5]\O{\Ribbon}\\
  \vstr{10}\step[.25]\id\step[1.5]\id\\[-2pt]
  \obox{2}{\mathsf{B}^2}\\
  \vstr{33}\step[.5]\id\step[1]\id
  \end{tangles}\ \ =
  \begin{tangles}{l}
    \vstr{50}\step[.5]\id\step[1]\id\\
    \obox{2}{\Ribbon}\\ 
    \vstr{50}\step[.5]\id\step[1]\id
  \end{tangles}
\end{equation}
Whenever it exists, choosing $\vartheta=\Ribbon$ in~\eqref{YD-loop}
makes $\chi_{\RepZ}$ ``multiplicative'' in $\RepZ$.  To show this, we
calculate $\chi_{\RepW}(\chi_{\RepZ}(\RepY))$ by sliding one of the
diagrams along the $\RepY$ line into the middle of the other and then
expanding:
\begin{equation}\label{B--B}
  \kern-6pt
  \begin{tangles}{l}
    \id\fobject{\kern10pt\RepY}\step[3]\fobject{\kern-12pt\RepZ}\coev
    \fobject{\kern12pt\DualZ}\\[-3.9pt]
    \vstr{50}\dh\step[1.5]\dd\step[2]\id\\[-2pt]
    \step[.25]\obox{2}{\mathsf{B}^2}\step[2.75]\id\\[-2pt]
    \vstr{50}\hdd\step[1.5]\d\step[2]\id\\
    \id\step[3]\O{\theta}\step[1]\dd\\
    \vstr{50}\id\step[3]\hx\\
    \id\step[3]{\makeatletter\@ev{0,\hm@de}{5,\hm@detens}{10}b\makeatother}\\
    \id\fobject{\kern10pt\RepY}\step[3]\fobject{\kern-14pt\RepW}\coev
    \fobject{\kern14pt\DualW}\\[-.1pt]
    \vstr{50}\dh\step[1.5]\dd\step[2]\id\\[-2pt]
    \step[.25]\obox{2}{\mathsf{B}^2}\step[2.75]\id\\[-2pt]
    \vstr{50}\hdd\step[1.5]\d\step[2]\id\\
    \id\step[3]\O{\theta}\step[1]\dd\\
    \vstr{50}\id\step[3]\hx\\
    \id\step[3]{\makeatletter\@ev{0,\hm@de}{5,\hm@detens}{10}b\makeatother}
  \end{tangles}
  \qquad = \ \
  \begin{tangles}{l}
    \step[.5]\ld[2]\step[2]    
    \\
    \hcd\step[1.5]\x
    {\makeatletter\@ev{0,\hm@de}{25,\hm@detens}{50}t\makeatother}
    \step[5]\id
    \\
    \id\step\id\step[1]\hld\step[1]\ld
    \step[5]\id
    \step[-2]\hh\coev\\
    \id\step\id\step[1]\id\step[.5]\id\step[1]\id\step\d
    \step[1]\dd\step[1]\id\step[1]\id\\
    \id\step\hx\step[.5]\id\step[1]\O{\medA}
    \step[1.5]\obox{2}{\mathsf{B}^2}\step[1.5]\id\step[1]\id
    \\[-.8pt]
    \id\step\id\step[1]\id\step[.5]\id\step[1]\id\step\dd
    \step[1]\hd\step[1.5]\id\step[1]\id\\[-.9pt]
    \hcu\step[1]\hlu\object{\raisebox{13.1pt}{\kern-4pt\tiny$\blacktriangleright$}}\step\lu\object{\raisebox{8pt}{\kern-4pt\tiny$\blacktriangleright$}}
    \step[2.5]\O{\theta}\step[1]\ddh\step[1]\id\\
    \step[.5]\id\step[2]\x
    \step[2.5]\hx\step[1.5]\id\\
    \step[.5]\lu[2]\object{\raisebox{8pt}{\kern-4pt\tiny$\blacktriangleright$}}\step[2]\O{\theta}
    \step[2.5]{\makeatletter\@ev{0,\hm@de}{5,\hm@detens}{10}b\makeatother}
    \step[1.5]\ne{3}\\
    \step[2.5]\id\step[2]\x\\
    \step[2.5]\id\step[2]\ev
  \end{tangles}
  \ \ \ = \quad
  \begin{tangles}{l}
    \vstr{90}\step[1]\ld[2]\step[2]\Coev\\
    \vstr{90}\step[.5]\hcd\step[1.5]\x
    \step[3]\d\\
    \vstr{90}\ddh\step[1]\id\step[1]\hld\fobject{\checkmark}\step[2]\d\step[3]\d
    \step[-2]
    \hh\coev\\
    \vstr{90}\id\step[1.5]\hx\step[.5]\dh\step[2.5]\x\step[1]\id\step\id\\
    \id\step[1.5]\id\step[.5]\hcd\step[.5]\hd\step[2]\id
    \step\ld\step[1]\id\step\id\\
    \vstr{90}\id\step[1.5]\id\step[.5]\id\step[1]\hx\step\ld\fobject{\kern-6pt\checkmark}
    \step\O{\medA}\step\id\step[1]\id\step\id\\
    \vstr{90}\id\step[1.5]\id\step[.5]\lu\object{\raisebox{7pt}{\kern-4pt\tiny$\blacktriangleright$}}\step\hx\step\id\step\lu\object{\raisebox{7pt}{\kern-4pt\tiny$\blacktriangleright$}}\step[1]\id\step\id
    \\
    \vstr{90}\dh\step[1]\hd\step[1]\hx\step[1]\lu\object{\raisebox{7pt}{\kern-4pt\tiny$\blacktriangleright$}}\step[2]\id\step[1]\id\step\id\\
    \vstr{90}\step[.5]\d\step[.5]\hcu\step\id\step[2]\x\step[1]\id\step\id\\
    \vstr{90}\step[1.5]\hcu\step[1.5]\x\step[2]\O{\theta}\step\id\step\id\\
    \vstr{90}\step[2]\lu[2]\object{\raisebox{7pt}{\kern-4pt\tiny$\blacktriangleright$}}\step[2]\O{\theta}\step[2]\hx\step[1]\id\\
    \vstr{90}\step[4]\id\step[2]\d\step{\makeatletter\@ev{0,\hm@de}{5,\hm@detens}{10}b\makeatother}\step\dd\\
    \vstr{90}\step[4]\id\step[3]\x\\
    \vstr{90}\step[4]\step[3]\ev
  \end{tangles}
\end{equation}
In the last diagram, we recognize the diagonal coaction (the two
$\checkmark$) and action (two $\blacktriangleright$ just below the
respective checkmarks) on a tensor product of two Yetter--Drinfeld
modules, as in~\eqref{act-coact-YD}.  In the bottom right part of the
diagram, we recall that \quad $\begin{tangles}{l} \O{\theta}
\end{tangles}\ \ = \ \ \begin{tangles}{l}
  \O{\Ribbon}\\
  \O{\Times}
\end{tangles}$\quad and calculate
\begin{equation*}
  \begin{tangles}{l}
    \step[-.5]\O{\Times}\step[1.5]\O{\Times}\step[1.5]\id\step\id\\
    \step[-.5]\hd\step\id\step[1.5]\id\step\id\\
    \id\step\id\step\ddh\step[1]\id\\
    \id\step\hx\step\ddh\\
    \d{\makeatletter\@ev{0,\hm@de}{5,\hm@detens}{10}b\makeatother}
    \step[1]\dd\\
    \step\hx\\
    \step{\makeatletter\@ev{0,\hm@de}{5,\hm@detens}{10}b\makeatother}
  \end{tangles}
  \ \ = \ \
  \begin{tangles}{l}
    \step[-.5]\O{\Times}\step[1.5]\O{\Times}\step[1.5]\id\step\id\\
    \step[-.5]\hd\step\id\step[1]\hdd\step[.5]\ddh\\
    \hx\step[1]\hxx\\
    \id\step\hx\step\id\\
    \d{\makeatletter\@ev{0,\hm@de}{5,\hm@detens}{10}b\makeatother}
    \step[1]\dd\\
    \step[1]\hx\\
    \step{\makeatletter\@ev{0,\hm@de}{5,\hm@detens}{10}b\makeatother}
  \end{tangles}
  \ \ = \ \
  \begin{tangles}{l}
    \step[-.5]\O{\Times}\step[1.5]\O{\Times}\step[1.5]\id\step\id\\
    \step[-.5]\hd\step\id\step[1]\hdd\step[.5]\ddh\\
    \hx\step[1]\hxx\\
    \hx\step[1]\id\step[1]\id\\
    \id\step[1]\hx\step[1]\id\\
    \hx\step[1]\hx\\
    {\makeatletter\@ev{0,\hm@de}{5,\hm@detens}{10}b\makeatother}
    \step[2]{\makeatletter\@ev{0,\hm@de}{5,\hm@detens}{10}b\makeatother}
  \end{tangles}
  \ \ = \ \
  \begin{tangles}{l}
    \step[-.5]\O{\Times}\step[1.5]\O{\Times}\step[1.5]\id\step\id\\
    \step[-.5]\hd\step\id\step[1]\hdd\step[.5]\ddh\\
    \hx\step\id\step\id\\
    \hx\step\id\step\id\\
    \id\step\hx\step\id\\
    \hx\step\hx\\
    \id\step\hx\step\id\\
    \Ev{\makeatletter\@ev{0,\hm@de}{5,\hm@detens}{10}b\makeatother}
  \end{tangles}
  \ \ 
  = \ \
  \begin{tangles}{l}
    \vstr{67}\step[.5]\id\step[1]\id\step[1]\id\step[1]\id\\
    \vstr{67}\step[.5]\hx\step[1]\id\step[1]\id\\[-1pt]
    \hld\step[1]\id\step[1]\id\step[1]\id\\
    \vstr{67}\id\step[.5]\hx\step[1]\id\step[1]\id\\
    \vstr{67}\id\step[.5]\id\step[1]\hx\step[1]\id\\
    \vstr{67}\id\step[.5]\hx\step[.5]\step[.5]\hx\\[-1pt]
    \hlu\object{\raisebox{13.1pt}{\kern-4pt\tiny$\blacktriangleright$}}
    \step[.5]\hld\step[1]\id\step[1]\id\\
    \vstr{67}\step[.5]\id\step[.5]\id\step[.5]\hx\step[1]\id\\[-1pt]
    \step[.5]\id\step[.5]\hlu\object{\raisebox{13.1pt}{\kern-4pt\tiny$\blacktriangleright$}}\step[1]\id\step[1]\id\\
    \step[.5]\Ev{\makeatletter\@ev{0,\hm@de}{5,\hm@detens}{10}b\makeatother}
  \end{tangles}
  \ \ \stackrel{\checkmark}{=} \ \ 
  \begin{tangles}{l}
    \vstr{90}\hld\step[1]\id\step[1]\id\step[1]\id\\
    \vstr{90}\vstr{50}\id\step[.5]\hx\step[1]\id\step[1]\id\\[-1pt]
    \vstr{90}\hlu\object{\raisebox{11.5pt}{\kern-4pt\tiny$\blacktriangleright$}}
    \step[-.5]\hld\step[1]\id\step[1]\id\step[1]\id\\
    \vstr{50}\id\step[.5]\hx\step[1]\id\step[1]\id\\[-1pt]
    \vstr{90}\hlu\object{\raisebox{11.5pt}{\kern-4pt\tiny$\blacktriangleright$}}
    \step[.5]\hld\step[1]\id\step[1]\id\\
    \vstr{50}\step[.5]\id\step[.5]\id\step[.5]\hx\step[1]\id\\[-1pt]
    \hld\step[.5]\hlu\object{\raisebox{13.1pt}{\kern-4pt\tiny$\blacktriangleright$}}\step[.5]\hld\step[1]\id\\
    \vstr{50}\id\step[.5]\hx\step[.5]\id\step[.5]\hx\\[-1pt]
    \vstr{90}\hlu\object{\raisebox{11.5pt}{\kern-4pt\tiny$\blacktriangleright$}}
    \step[.5]\hld\step[.5]\hlu\object{\raisebox{11.5pt}{\kern-4pt\tiny$\blacktriangleright$}}\step[1]\id\\
    \vstr{50}\step[.5]\id\step[.5]\id\step[.5]\hx\step[1]\id\\[-1pt]
    \vstr{90}\step[.5]\id\step[.5]\hlu\object{\raisebox{11.5pt}{\kern-4pt\tiny$\blacktriangleright$}}\step[1]\id\step[1]\id\\
    \vstr{90}\step[.5]\Ev{\makeatletter\@ev{0,\hm@de}{5,\hm@detens}{10}b\makeatother}
  \end{tangles}
  \ \ = \ \
  \begin{tangles}{l}
    \step[.5]\id\step[1]\id\step[1]\id\step[1]\id\\
    \obox{2}{\mathsf{B}^2}\step[.5]\id\step[1]\id\\
    \step[.5]\id\step[1]\id\step[1]\id\step[1]\id\\
    \obox{2}{\TIMES}\step[.5]\id\step[1]\id\\
    \step[.5]\id\step\hx\step\id\\
    \step[.5]\hx\step\hx\\
    \step[.5]\id\step\hx\step\id\\
    \step[.5]\Ev{\makeatletter\@ev{0,\hm@de}{5,\hm@detens}{10}b\makeatother}
  \end{tangles}
\end{equation*}

\noindent
where the first three equalities are elementary (and well-known)
rearrangements, the fourth involves~\eqref{equiv-to-2}, and the
checked equality is verified by repeatedly applying the
Yetter--Drinfeld axiom in its right-hand side.  The sixth diagram
involves $\mathsf{B}^2$ in the upper part and the diagonal action and
coaction~\eqref{act-coact-YD} in the lower part, which gives the last
equality.  We therefore conclude that if~\eqref{Rib-def} holds,
then
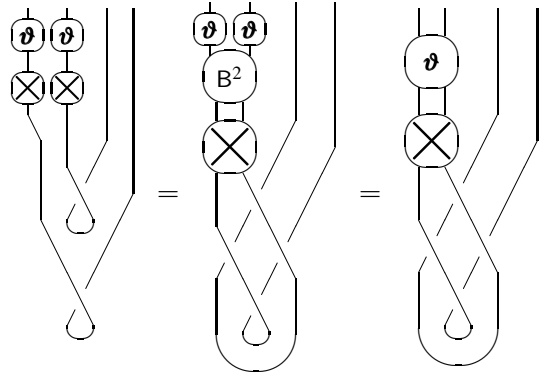
\begin{equation*}
  \begin{tangles}{l}
    \step[-.5]\O{\Ribbon}\step[1.5]\O{\Ribbon}\step[1.5]\id\step\id\\
    \step[-.5]\O{\Times}\step[1.5]\O{\Times}\step[1.5]\id\step\id\\
    \step[-.5]\hd\step\id\step[1]\ddh\step\id\\
    \id\step\hx\step\ddh\\
    \d{\makeatletter\@ev{0,\hm@de}{5,\hm@detens}{10}b\makeatother}
    \step[1]\dd\\
    \step\hx\\
    \step{\makeatletter\@ev{0,\hm@de}{5,\hm@detens}{10}b\makeatother}
  \end{tangles}
  \ \ = \ \
  \begin{tangles}{l}
    \step[.25]\O{\Ribbon}\step[1.5]\O{\Ribbon}\step[1.75]\id\step[1.5]\id\\[-2pt]
    \obox{2}{\mathsf{B}^2}\step[1.5]\id\step[1.5]\id\\
    \vstr{33}\step[.5]\id\step[1]\id\step[2]\id\step[1.5]\id\\
    \obox{2}{\TIMES}\step[.5]\dd\step[1]\ddh\\
    \step[.5]\id\step\hx\step\dd\\
    \step[.5]\hx\step\hx\\
    \step[.5]\id\step\hx\step\id\\
    \step[.5]\Ev{\makeatletter\@ev{0,\hm@de}{5,\hm@detens}{10}b\makeatother}
  \end{tangles}
  \ \ = \ \
  \begin{tangles}{l}
    \vstr{50}\step[.5]\id\step[1]\id\step[2]\id\step[1.5]\id\\
    \obox{2}{\Ribbon}\step[1.5]\id\step[1.5]\id\\
    \vstr{50}\step[.5]\id\step[1]\id\step[2]\id\step[1.5]\id\\
    \obox{2}{\TIMES}\step[.5]\dd\step[1]\ddh\\
    \step[.5]\id\step\hx\step\dd\\
    \step[.5]\hx\step\hx\\
    \step[.5]\id\step\hx\step\id\\
    \step[.5]\Ev{\makeatletter\@ev{0,\hm@de}{5,\hm@detens}{10}b\makeatother}
  \end{tangles}
\end{equation*}

\noindent
Substituting this in~\eqref{B--B} shows that $\chi$ is indeed
``multiplicative'':
$\chi_{\RepW}(\chi_{\RepZ}(\RepY))=\chi_{\RepW\tensor\RepZ}(\RepY)$.

\section{Rank-one Nichols algebra}\label{sec:p}
We specialize the preceding sections to the case of a rank-one Nichols
algebra $\Nich\!_p$, whose relation to the $(p,1)$ logarithmic CFT{}
models was emphasized in~\cite{[STbr]}.  An integer $p\geq 2$ is fixed
throughout.

\subsection*{Notation}
We fix the primitive $2p$th root of unity
 \begin{equation*}
  \q=e^{\frac{i\pi}{p}}
\end{equation*}
and introduce the $q$-binomial coefficients
\begin{equation*}
  \Abin{r}{s}=\ffrac{\Afac{r}}{\Afac{s}\Afac{r-s}},
  \quad\Afac{r}=\Aint{1}\dots\Aint{r},
  \quad\Aint{r}=\ffrac{q^{2r}-1}{q^2-1},
\end{equation*}
which are assumed to be specialized to $q=\q$.

We sometimes use the notation $(a)_N=
a\;\mathrm{mod}\;N\in\{0,1,\dots,N-1\}$.

\subsection{The braided Hopf algebra
  $\Nich\!_p$}\label{sec:screenings}
The rank-$1$ Nichols algebra $\Nich\!_p$ is $\Nich(X)$ for a
one-dimensional braided linear space $X$.  We fix an element $F$ (a
single screening in the CFT{} language) as a basis in $X$.  The
braiding, taken from CFT{}, is
\begin{equation}\label{Psi-F-F}
  \Psi(F(r)\tensor F(s))=\q^{2 r s} F(s)\tensor F(r).
\end{equation}
Shuffle product~\eqref{Sh-prod} then becomes
\begin{equation*}
  F(r)\,F(s)=\Abin{r+s}{r}F(r+s)
\end{equation*}
and coproduct~\eqref{cutting} is
$\Delta:F(r)\mapsto\sum\limits_{s=0}^r F(s)\tensor F(r-s)$.  The
antipode defined in~\eqref{antipode} acts as $\A(F(r))=(-1)^r
\q^{r(r-1)}F(r)$.

The algebra $\Nich\!_p$ is the linear span of $F(r)$ with $0\leq r\leq
p-1$.  It can also be viewed as generated by a single element $F$,
such that $F(r)=\ffrac{1}{\Afac{r}}F^r$, $r\leq p-1$, with $F^p=0$.
We write $F=F(1)$.

Because $X$ is now one-dimensional, we can think of
$\xymatrix@1@C=20pt{ \ar@{--}|(.55){\cross}[0,2]&& }$ as just $F$, and
write
\begin{equation*}
  F(r)={}\ \xymatrix@1@C=40pt{
    \ar@{--}|(.2){\cross}|(.55){\cross}|(.8){\cross}[0,2]&& }
  \qquad\text{($r$ crosses)}.
\end{equation*}

\subsection{Yetter--Drinfeld $\Nich\!_p$-modules}
We specialize the construction of Yetter--Drinfeld $\Nich(X)$-modules
in Sec.~\ref{sec:Nich} to $\Nich\!_p$.  The construction involves
another braided vector space $Y$, a linear span of vertex operators
present in the relevant CFT{} model.

\label{sec:dressedVO}
\subsubsection{The vertices}
For the $(p,1)$ model corresponding to $\Nich\!_p$
(see~\cite{[FHST]}), $Y$ is a $2p$-dimensional space
\begin{equation*}
  Y= \text{span}(V^a\mid a\in\oZ_{4p})
\end{equation*}
with the diagonal braiding
\begin{equation}\label{Psi-V-V}
  \Psi(V^a\tensor V^{b}) = \q^{\frac{a b}{2}}  V^{b}\tensor V^a.
\end{equation}
and with
\begin{equation}\label{Psi-F-V}
  \Psi(V^a\tensor\thev{r})=\q^{-a r}\thev{r}\tensor V^a,
  \quad
  \Psi(\thev{r}\tensor V^a)=\q^{-a r}\,V^a\tensor\thev{r}.
\end{equation}
This suffices for calculating the ``cumulative adjoint'' $\Nich\!_p$
action on multivertex Yetter--Drinfeld modules, as we describe next.

In what follows, the integers $a$, $b$, \dots\ are tacitly considered
modulo $4p$.

\subsubsection{Multivertex Yetter--Drinfeld $\Nich\!_p$-modules}
We saw in Sec.~\ref{sec:Nich} that multivertex Yet\-ter--Drinfeld
modules (see~\eqref{ex1} and~\eqref{2items}) can be represented as an
essentially ``combinatorial'' construction for the crosses to
populate, in accordance with the braiding rules, line segments that
are separated from one another by vertex operators, e.g.,
$\xymatrix@R=6pt@C=50pt@1{ \ar@*{[|(1.6)]}@{-}|(.15){\cross}
  |(.25){\punct}|(.4){\cross}|(.5){\cross}|(.65){\punct} |(.8){\cross}
  |(.90){\punct}[0,2]&& }$, where ${\times}=X$ and ${\circ}=Y$ (for a
finite-dimensional Nichols algebra, each ``segment'' can carry only
finitely many crosses).  In the rank-$1$ case, each cross can be
considered to represent the $F$ element, and each segment is fully
described just by the number of the $F$s sitting there.  For example,
each two-vertex Yetter--Drinfeld module is a linear span of
\begin{equation}\label{Vsatb}
   \VertexIII{s}{a}{t}{b}
   =\xymatrix@R=4pt@C=70pt{
    \ar@*{[|(1.6)]}@{-}|(.1){\cross}|(.25){\cross}|(.5){\punct}^(.5){{}^{\scriptstyle
        a}}|(.65){\cross}|(.8){\punct}^(.8){{}^{\scriptstyle
        b}}[0,2]&&}
  \!\!\!\!,
\end{equation}

\noindent
where $s$ and $t$ must not exceed $p-1$ ($s=2$ and $t=1$ in the
picture) and $a$ and $b$ indicate $V^a$ and $V^b$.  Because the
braiding is diagonal, there is a $\Nich\!_p$ module comodule for each
fixed $a$ and $b$ (and $c$, \dots\ for multivertex modules).

The simplest, one-vertex Yetter--Drinfeld $\Nich\!_p$-modules are
spanned by
\begin{equation}\label{Vsa}
  \VertexI{s}{a}=\xymatrix@R=4pt@C=70pt{
    \ar@*{[|(1.6)]}@{-}|(.1){\cross}|(.35){\cross}|(.55){\cross}|(.7){\punct}^(.7){{}^{\scriptstyle a}}[0,2]&&\text{\qquad\small($s$ crosses),}}
\end{equation}

\noindent
where $s$ 
ranges over a subset of $[0,\dots,p-1]$.  The $\Nich\!_p$ coaction is by
``deconcatenation up to the first vertex'' in all cases, i.e.,
\begin{align*}
  \delta\,\VertexI{s}{a} &= \sum_{r=0}^s
  F(r)
  \tensor\VertexI{s-r}{a},
  \\
  \delta \VertexIII{s}{a}{t}{b}
  &=\sum_{r=0}^s
  F(r)
  \tensor\VertexIII{s-r}{a}{t}{b},
\end{align*}
and similarly for $\VertexV{s}{a}{t}{b}{u}{c}$, and so on.

The $\Nich\!_p$ action (which is the left adjoint action~\eqref{adja})
is then calculated as
\begin{align*}
  F\adjoint
  \VertexI{s}{a} &=
  \xi
  \Aint{s - a}
  \Aint{s + 1} 
  \VertexI{s + 1}{a},\qquad \xi = 1-\q^2,
  \\
  \intertext{and the cumulative adjoint evaluates on multivertex
    spaces as}
  F\adjoint
  \VertexIII{s}{a}{t}{b} &= 
  \xi
  \Aint{s + 2 t - a - b}
  \Aint{s + 1}
  \VertexIII{s + 1}{a}{t}{b}
  +
  \xi
  \q^{2 s - a}
  \Aint{t - b}
  \Aint{t + 1}
  \VertexIII{s}{a}{t + 1}{b},
  \\
  F\adjoint \VertexV{s}{a}{t}{b}{u}{c}
  &=
  \xi
  \Aint{s + 2 t + 2 u - a - b - c}
  \Aint{s + 1} \VertexV{s + 1}{a}{t}{b}{u}{c}
  \\
  &{}+ \q^{2 s - a}
  \xi
  \Aint{t + 2 u - b - c}
  \Aint{t + 1}
  \VertexV{s}{a}{t + 1}{b}{u}{c}
  + \q^{2 s + 2 t - a - b}
  \xi
  \Aint{u - c}
  \Aint{u + 1}
  \VertexV{s}{a}{t}{b}{u + 1}{c},
\end{align*}
and so on.

The braiding follows from~\eqref{Psi-F-F}, \eqref{Psi-V-V},
and~\eqref{Psi-F-V}, for example,
\begin{equation}\label{braiding11}
  \Psi(\VertexI{s}{a}\tensor \VertexI{t}{b}) =
  \q^{\half(a - 2 s) (b - 2 t)}
  \VertexI{t}{b}\tensor\VertexI{s}{a}.
\end{equation}

\subsection{Module types and decomposition}
We now study the category of Yetter--Drinfeld $\Nich\!_p$-modules in
some detail: we find how the one-vertex and two-vertex spaces
decompose into indecomposable Yetter--Drinfeld $\Nich\!_p$-modules.  We
first forget about braiding and study only \textit{the module comodule
  structure}; the action and coaction \textit{are} related by the
Yetter--Drinfeld axiom, but we try to avoid speaking of
Yetter--Drinfeld modules before we come to the braiding.

\subsubsection{}The relevant module comodules, which we construct
explicitly in Appendix~\bref{app:modules}, are as follows:
\begin{itemize}
\item simple $r$-dimensional module comodules $\RepX(r)$, $1\leq r\leq
  p$; for $r=p$, we sometimes use the special notation
  $\RepS(p)=\RepX(p)$;

\item the $p$-dimensional extensions
  \begin{align}\label{V-def}
    \raisebox{-.7\baselineskip}{$\RepV[r]\;={}$} &\xymatrix@=15pt{
      &\RepX(p-r)\ar@/_12pt/[dl]\\
      \RepX(r)& }
    \qquad\qquad\raisebox{-.7\baselineskip}{$1\leq r\leq p-1,$}
  \end{align}
  where the arrow means that $\delta\RepX(p-r)\subset($the ``trivial''
  piece~$\Nich\!_p\tensor \RepX(p-r))+\Nich\!_p\tensor
  \RepX(r)$.

\item $2p$-dimensional indecomposable module comodules $\RepP[r]$ with
  the structure of subquotients
  \begin{align}
    \label{P-def}
    \raisebox{-1.55\baselineskip}{$\RepP[r]\;={}$}
    &\xymatrix@=12pt{
      &\modX(p-r)
      \ar@/_12pt/[dl]
      \ar[dr]
      \\
      \modX(r)\ar[dr] &&\modX(r)\ar@/_12pt/[dl]
      \\
      &\modX(p-r)&
    }\qquad\qquad\raisebox{-1.55\baselineskip}{$1\leq r\leq p-1.$}
  \end{align}
\end{itemize}

\subsubsection{} 
We also show in Appendix~\bref{app:modules} that the $p^2$-dimensional
\textit{one-vertex space}
\begin{equation*}
  \mathbb{V}_p(1)\equiv
  \text{Span}(\VertexI{s}{a} \mid 0\leq a,s\leq p-1)
\end{equation*}
decomposes into $\Nich\!_p$ module comodules as
\begin{equation}\label{V1-decomp}
  \mathbb{V}_p(1) =
  \RepS(p)\;\oplus\bigoplus_{1\leq r\leq p-1}\RepV[r]
\end{equation}
and the $p^4$-dimensional \textit{two-vertex space}
\begin{equation*}
  \mathbb{V}_p(2)\equiv
  \text{Span}(\VertexIII{s}{a}{t}{b} \mid 0\leq a,b,s,t\leq p-1)
\end{equation*}
decomposes as
\begin{equation}\label{V2-decomp}
  \mathbb{V}_p(2) =
  p^2 \RepS(p)\;\oplus
  \bigoplus_{1\leq r\leq p-1} 2 r (p - r) \RepV[r]\;\oplus
  \bigoplus_{1\leq r\leq p-1}(p - r)^2 \RepP[r].
\end{equation}
\textit{Multivertex spaces} give rise to ``zigzag'' Yetter--Drinfeld
modules, which we do not consider here.

\subsubsection{Notation}
Compared with representation theory of Lie algebras, the role of
highest-weight vectors is here played by \textit{left coinvariants}
$\VertexI{0}{a}$ and $\VertexIII{0}{a}{t}{b}$.  When a module comodule
of one of the above types $\RepA=\RepX$, $\RepV$, or $\RepP$ is
constructed starting with a left coinvariant, we use the notation
$\RepA_{0}^{\{a\}}$ or $\RepA_{0,\,t}^{\{a,\,b\}}$ to indicate the
coinvariant, and sometimes also use the notation such as
$\RepXii{0}{a}{t}{b}(r)$ to indicate the dimension (although it is
uniquely defined by $a$, $t$, $b$, and the module type).

\subsubsection{}\label{1-vertex}The module comodules that can be
constructed starting with one-vertex coinvariants $\VertexI{0}{a}$ are
classified immediately, as we show
in~\bref{app:1-vertex}.\pagebreak[3] The module comodule
\textit{generated} from $\VertexI{0}{a}$ under the $\Nich\!_p$ action
is isomorphic to $\RepX(r)$ whenever $(a)_p = r - 1$ ($1\leq r\leq
p$).  If $r\leq p-1$, then extension~\eqref{V-def} follows
immediately.

\subsubsection{}\label{sec:all-cases}The strategy to classify
two-vertex $\Nich\!_p$ module comodules according to their
characteristic left coinvariant $\VertexIII{0}{a}{t}{b}$ is to
consider the following cases that can occur under the action of $F(s)$
on the left coinvariant.
\begin{enumerate}
\item $F(s)\adjoint \VertexIII{0}{a}{t}{b}$ is nonvanishing and not a
  coinvariant for all $s$, $1\leq s\leq p-1$.  In this case, there are
  the possibilities that
  \begin{enumerate}
  \item\label{caseL} $F(s)\adjoint \VertexIII{0}{a}{t}{b}$ is a
    coinvariant, i.e., $F(s)\adjoint
    \VertexIII{0}{a}{t}{b}=\mathrm{const}\;\VertexIII{0}{a}{t+s}{b}$,
    for some $s\leq p-1$, and

  \item\label{caseS} $F(s)\adjoint \VertexIII{0}{a}{t}{b}$ is not a
    coinvariant for any $s\leq p-1$.
  \end{enumerate}

\item $F(s)\adjoint \VertexIII{0}{a}{t}{b}=0$ for some $s\leq p-1$.
  In this case, further possibilities are
  \begin{enumerate}
  \item\label{casenone} For some $s'<s$, \ $F(s')\adjoint
    \VertexIII{0}{a}{t}{b}$ is a coinvariant, and

  \item $F(s')\adjoint \VertexIII{0}{a}{t}{b}$ is not a coinvariant
    for any $s'<s$.  We then distinguish the cases where
    \begin{enumerate}
    \item\label{caseB} $\VertexIII{0}{a}{t}{b}$ is in the image of
      $F$, and
    \item\label{caseX} $\VertexIII{0}{a}{t}{b}$ is not in the image of
      $F$.
    \end{enumerate}
  \end{enumerate}
\end{enumerate}

We show in Appendix~\bref{app:modules} that these cases are resolved
as follows in terms of the parameters $a$, $t$, and $r=(a + b - 2 t)_p
+ 1$:
\begin{description}
\item[\normalfont\ref{caseL}] $1\leq r\leq p-1$ and either $t \leq
  (a)_p - r$ \ or \ $(a)_p + 1 \leq t \leq p - r - 1$. \ Then the left
  coinvariant is the leftmost coinvariant in~\eqref{P-mod-ind}, and
  the Yetter--Drinfeld module generated from it is the ``left--bottom
  half'' $\RepL(r)$ of $\RepP[r]$ (see~\bref{sec:LB}).

\item[\normalfont\ref{caseS}] $r=p$. \ Then $\RepX(p)\equiv\RepS(p)$
  is generated from the left coinvariant.

\item[\normalfont\ref{casenone}] is not realized.

\item[\normalfont\ref{caseB}] $1\leq r\leq p-1$ and either $t \geq p -
  r + (a)_p + 1$ \ or \ $p - r \leq t \leq (a)_p$. \ Then the bottom
  Yetter--Drinfeld submodule $\RepB(r)$ in $\RepP[p-r]$ is generated
  from the left coinvariant.

\item[\normalfont\ref{caseX}] $1\leq r\leq p-1$ and either $(t\leq
  (a)_p$ and $(a)_p - r + 1\leq t\leq p - 1 - r)$ \ or \ $(t\geq (a)_p
  + 1$ and $p - r\leq t\leq p - r + (a)_p)$. \ Then $\RepX(r)$ is
  generated from the left coinvariant.
\end{description}

\subsection{Braiding sectors}\label{sec:braiding-sectors}
The $\RepX(r)$ and the other module comodules appearing above satisfy
the Yetter--Drinfeld axiom.  Considering them as Yetter--Drinfeld
$\Nich\!_p$-modules means that isomorphic module comodules may be
distinguished by the braiding.
This is indeed the case: for example, shifting $a\to a+p$ in
$\RepXi{a}$ or $\RepXii{0}{a}{t}{b}$ does not affect the module
comodule structure described in Appendix~\bref{app:modules}, but
changes the braiding with elements of $\Nich\!_p$ by a sign in
accordance with~\eqref{Psi-F-V}.  We thus have \textit{pairs}
$(\RepA_\nu,\RepA_{\nu+1})$, $\nu\in\oZ_2$, of isomorphic module
comodules distinguished by a sign occurring in their braiding.  In
particular, there are $2p$ nonisomorphic simple Yetter--Drinfeld
modules.

Further, these Yetter--Drinfeld modules can be viewed as elements of a
braided category, whose braiding (see~\eqref{B+B2})
involves~\eqref{braiding11}.  The dependence on $a$
in~\eqref{braiding11} is modulo $4p$, and hence we have not pairs but
quadruples $(\RepA_\nu)_{\nu\in\oZ_4}$, with the different $\RepA_\nu$
distinguished by their braiding with other such modules.  In
particular, there are $4p$ nonisomorphic simple objects in this
braided category of Yetter--Drinfeld $\Nich\!_p$-modules~\cite{[STbr]}.

It is convenient to write $a=(a)_p-\nu p,$
$\nu\in\oZ_4$~\cite{[STbr]}, and introduce the notation
$\RepX(r)_{\nu}$ for simple modules,
with
\begin{equation*}
  \RepXi{a}\cong\RepX(r)_{\nu} \quad\text{whenever}\quad a=r-1-\nu p.
\end{equation*}
As before, $r$ is the dimension, and we sometimes refer to $\nu$ as
the braiding sector or braiding index.  For $\nu\in\oZ_4$, the
isomorphisms are in the braided category of Yetter--Drinfeld
$\Nich\!_p$-modules.  The ``quadruple structure'' occurs totally
similarly for other modules, including those realized in multivertex
spaces; for example, for any $a,b\in\oZ$, we have the isomorphisms
among the simple Yetter--Drinfeld modules realized in the two-vertex
space (cf.~\eqref{2X-mod}!):
\begin{equation*}
  \RepXii{0}{a}t{b}\cong\RepX(r)_\nu
  \quad\text{whenever \ $a + b - 2 t = r-1 - \nu p$ \ and \
\eqref{cond:Xi}${}\lor{}$\eqref{cond:Xii} \ holds.}
\end{equation*}
For the reducible extensions as in~\eqref{V-def}, the two subquotients
have adjacent braiding indices, and we conventionally use one of them
in the notation for the reducible module:
\begin{equation}\label{X-X-nu}
  \raisebox{-.8\baselineskip}{\mbox{$\RepV[r]_{\nu}\
      =\ $}}
  \xymatrix@=15pt{
    &\RepX(p-r)_{\nu + 1}\ar@/_12pt/[dl]\\
    \RepX(r)_{\nu}&
    },
\end{equation}
and $\RepVii{0}{a}{t}{b}[r]_{\nu}\cong \RepV[r]_{\nu}$ whenever
$a+b-2t = r-1 -\nu p$ and \eqref{cond:Xi}${}\lor{}$\eqref{cond:Xii}
holds.

In~\eqref{P-def}, the relevant braiding indices range an interval of
three values, and we use the leftmost value in the notation for the
entire reducible Yetter--Drinfeld module, which
yields~\eqref{P-mod-ind}, with $\RepPii{0}{a}{t}{b}[r]_{\nu}
\cong{}\RepP[r]_{\nu}$ whenever $a+b-2t = r-1 -\nu p$ and
\eqref{ProjectiveQ} holds.

In the above formulas and diagrams, $\nu\in\oZ_4$ if the modules are
viewed as objects of the braided category of Yetter--Drinfeld
$\Nich\!_p$-modules.  But \textit{if the Yetter--Drinfeld
  $\Nich\!_p$-modules are considered as an entwined category, then the
  braiding sectors $\nu$ and $\nu+2$ become indistinguishable, and
  hence $\nu\in\oZ_2$}.  In particular, there are $2p$ nonisomorphic
simple objects in the entwined category of Yetter--Drinfeld
$\Nich\!_p$-modules.


\subsection{Proof of
  decomposition~\eqref{the-fusion}}Decomposition~\eqref{the-fusion}
can be derived from the list in~\bref{sec:all-cases} as follows.  The
fusion product~\eqref{fusion-def} of two one-vertex modules is the map
(assuming that $a,b\leq p-1$ to avoid writing $(a)_p$ and $(b)_p$)
\begin{equation}\label{VV-fusion}
  V^a_s\tensor V^b_t\mapsto
  \sum_{i=0}^bq^{-a i}\Abin{s+i}{s}\VertexIII{s+i}{a}{t-i}{b}.
\end{equation}
In evaluating $\RepXi{a}(s)\tensor\RepXi{b}(t)$, this formula is applied
for $0\leq s\leq a$ and $0\leq t\leq b$.  Then the left coinvariants
produced in the right-hand side are $\VertexIII{0}{a}{u}{b}$, where
$0\leq u\leq b$ \textit{and} $u\leq a$.  But the conditions defining
the different items in the list in~\bref{sec:all-cases} have the
remarkable property that the module $\RepA^{\{a,b\}}_{0,u}$ generated
from each such coinvariant is as follows:
\begin{equation}\label{the-property}
  \RepA^{\{a,b\}}_{0,u}=
  \begin{cases}
    \RepXii{0}{a}{u}{b},&a + b \leq p-1,\\
    \RepXii{0}{a}{u}{b},&a + b \geq p\ \text{ and } \ u \geq a + b
    - p + 2,\\
    \RepLii{0}{a}{u}{b},&a + b \geq p\ \text{ and } \ a + b - 2 u -
    p \geq 0,\\
    \RepSii{0}{a}{u}{b},&a + b \geq p\ \text{ and } \ a + b - 2 u -
    p = -1,\\
    \RepBii{0}{a}{u}{b},&a + b \geq p\ \text{ and } \ a + b - 2 u -
    p \leq -2.
  \end{cases}
\end{equation}
This is established (\textit{only} for $0\leq u\leq a,b\leq p-1$) by
direct inspection of each case in the list at the end
of~\bref{sec:all-cases}. \ The module $\RepLii{0}{a}{u}{b}$ is the
``left--bottom half'' of $\RepPii{0}{a}{u}{b}$, and
$\RepBii{0}{a}{u}{b}$ is the bottom sub(co)module in another $\RepP$
module; the details are given in~\bref{sec:LB}.  The crucial point is
that each $\RepLii{0}{a}{u}{b}$ can be extended to
$\RepPii{0}{a}{u}{b}$ (while the $\RepB$, on the other hand, are not
interesting in that they are sub(co)modules in the $\RepL$ that are
already present).  We next claim that each of the $\RepL$s occurring
in $\RepXi{a}(s)\tensor\RepXi{b}(t)$ indeed occurs there together with
the entire $\RepP$ module; this follows from counting the dimensions
and from the fact that there are no more left coinvariants among the
$\VertexIII{v}{a}{w}{b}$ appearing in the right-hand side
of~\eqref{VV-fusion} (and, of course, from the structure of the
modules described in Appendix~\bref{app:modules}).

Once it is established that each $\RepL$ occurs
in~\eqref{the-property} as a sub(co)module of the
corresponding~$\RepP$, it is immediate to see
that~\eqref{the-property} is equivalent to~\eqref{the-fusion}.

\subsection{Duality} We now recall Sec.~\ref{sec:duality}.  The
structures postulated there are indeed realized for the $n$-vertex
Yetter--Drinfeld $\Nich\!_p$-modules.

\subsubsection{One-vertex modules: $\mathrm{coev}$ and $\mathrm{ev}$
  maps}
\label{coev-1}
For the irreducible Yetter--Drinfeld module $\RepXi{a}\cong\RepX(r)$
as in~\bref{sec:Xi}, the coevaluation map
$\mathrm{coev}:k\to\RepXi{a}\tensor\dualXi{a}$ is given in terms of
dual bases as
\begin{equation*}
  \begin{tangles}{l}
    \vphantom{x}\\
     \Coev\\
    \fobject{\RepXi{a}}\fobject{\kern60pt\dualXi{a}}
  \end{tangles}\kern30pt
  \ = \sum_{s=0}^{r-1}\VertexI{s}{a}\tensor\UertexI{s}{-a},
  \qquad r = (a)_p+1,
\end{equation*}
and the evaluation map $\mathrm{ev}: \dualXi{a}\tensor\RepXi{a}\to k$,
accordingly, as
\begin{equation*}
  \begin{tangles}{l}
    \vphantom{x}\\
    \fobject{\dualXi{a}}\fobject{\kern60pt\RepXi{b}}\\[4pt]
    \Ev
  \end{tangles}\kern30pt
  \ : \
  \UertexI{s}{-a}\tensor\VertexI{t}{b}\mapsto 1
  \delta_{s,t}\delta_{a,b}.
\end{equation*}
We then use~\eqref{on-dual} to find the $\Nich\!_p$ module comodule
structure on the $U_s^a$.  Simple calculation shows that
\begin{align*}
  F(r)
  U_s^a &=
  \q^{r(r-1) - r a - 2 r s}
  \Abin{s}{r} (-\xi)^{r} \prod_{t=s-r}^{s-1}\Aint{t+a}\;U^a_{s-r},
  \\
  \delta U^a_s &= \sum_{r=0}^{p-1-s}
  (-1)^r\q^{-r a - 2 s r - r(r - 1)}
  F(r)\tensor U_{s+r}^a.
\end{align*}
It follows that we can identify $U_{s}^{a} =
(-1)^{a+s}\q^{(s+1)(s+a-2)} V_{p - 1 - s}^{a - 2 + 2p}$ (the action
and coaction\linebreak[0]---and in fact the braiding---are identical
for both sides).  The coevaluation and evaluation maps can therefore
be expressed as
\begin{equation*}
  \begin{tangles}{l}
    \vphantom{x}\\
    \Coev\\
    \fobject{\RepXi{a}}\step[3]\fobject{{}^\vee\!\RepXi{a}}
   \end{tangles}\kern15pt =
   \sum_{s=0}^{r-1}\VertexI{s}{a}\tensor\VertexI{p-1-s}{2p-a-2}
   (-1)^{a + s} \q^{(s+1) (s - a - 2)},\quad r=(a)_p+1,
\end{equation*}
and
\begin{equation*}
  \begin{tangles}{l}\Ev
  \end{tangles}\kern25pt
  :\ \VertexI{s}{a}\tensor\VertexI{t}{b}\mapsto
  \eval{\VertexI{s}{a}}{\VertexI{t}{b}}=
  (-1)^s\q^{-s^2+s(a-1)}
  \delta_{s+t,p-1}\delta_{a+b,2p-2}.
\end{equation*}

For $a\not\equiv p-1\;\text{mod}\;p$, evidently, $a = r-1 - \nu
p$ implies that $2 p - a - 2= p - r - 1 + (\nu + 1) p$, and
therefore the module left dual to $\RepV[r]$ in~\eqref{X-X-nu}, with
$r=(a)_p+1$, can be identified as
\begin{equation*}
  \raisebox{-\baselineskip}{\mbox{${}^{\vee}(\RepV[r]_{\nu})={}$}}
  \xymatrix@=15pt{
    .&\RepX(r)_{-\nu}\ar@/_12pt/[dl]\\
    \RepX(p-r)_{-\nu-1}&}
\end{equation*}
where $\RepX(r)_{-\nu}$ is dual to $\RepX(r)_{\nu}$ in~\eqref{X-X-nu}.

\textit{The properties expressed in~\eqref{equiv-to-2}
  and~\eqref{equiv-to-2} now hold}, as is immediate to verify.

\subsubsection{Two-vertex modules}\label{coev-2}
Similarly to~\bref{coev-1}, for the $\UertexIII{s}{a}{t}{b}$ that are
dual to the two-vertex basis,
\begin{equation*}
  \eval{\UertexIII{s}{a}{t}{b}}{\VertexIII{u}{c}{v}{d}}
  =\delta_{a+c,0}\delta_{b+d,0}\delta_{s,u}\delta_{t,v},
\end{equation*}
it follows from~\eqref{on-dual} that
\begin{equation*}
  F(r)\UertexIII{s}{a}{t}{b}
  =\sum_{u=0}^r(-1)^r \q^{r(r-1)-r(b+2s+2t)}
  \coeff{-a}{-b}{s-r+u}{t-u}{r}{u}\UertexIII{s-r+u}{a}{t-u}{b}.
\end{equation*}
Replacing here $s\to p-1-s$ and $t\to p-1-t$ and noting that the
coefficients $\coeff{a}{ b}{ s}{ t}{ r}{ u}$ in~\eqref{c-st} have the
symmetry
\begin{equation*}
  \coeff{a}{ b}{ s}{ t}{ r}{ u} =
 \q^{2 r (r + 2 t + 2 s - a - b)}
  \coeff{-a - 2}{ -b - 2}{ p - 1 - s - r + u}{ p - 1 - t - u}{ r}{ u},
  \quad r\geq u,
\end{equation*}
we arrive at the identification
\begin{equation*}
  \UertexM{a}{b}{s}{t} = (-1)^{t + s}
   \q^{(t + s + 2) (2 a + b + t + s - 3)}
   \VertexM{a - 2}{b - 2}{p - 1 - s}{p - 1  - t}.
\end{equation*}
Hence, under the pairing
\begin{equation*}
  \eval{\VertexIII{s}{a}{t}{b}}{\VertexIII{u}{c}{v}{d}}
  =(-1)^{s + t} \q^{(s + t) (2 a + b + 1 - s - t)}
  \delta_{a+c,-2}\delta_{b+d,-2}\delta_{s+u,p-1}\delta_{t+v,p-1},
\end{equation*}
the module left dual to $\RepPii{0}{a}{t}{b}$ can be identified with
$\RepPii{0}{-a-2}{p-r-t-1}{-b-2}$ (as before, $a+b-2t=r-1-\nu p$,
$1\leq r\leq p-1$).  The module dual to~\eqref{P-mod-ind} has the
structure
\begin{equation*}
    \raisebox{-2\baselineskip}{\mbox{%
        ${}^{\vee}\bigl(\RepPii{0}{a}{t}{b}[r]_{\nu}\bigr) =
        \RepPii{0}{-a-2}{p-r-t-1}{-b-2}[r]_{-2-\nu}\ \ =$}}
  \xymatrix@=12pt{
    &&\modX(p-r)_{-\nu-1}
    \ar@/_12pt/[dl]
    \ar[dr]
    &\\
    &\modX(r)_{-2-\nu}\ar[dr] &&\modX(r)_{-\nu}\ar@/_12pt/[dl]
    \\
    &&\modX(p-r)_{-\nu-1}&
  }
\end{equation*}

\subsection{Ribbon structure}
We set
\begin{equation*}
  \Ribbon\, V^{a}_s = \q^{\half((a + 1)^2 - 1)}V^{a}_s,
\end{equation*}
which obviously commutes with the $\Nich\!_p$ action and coaction, and
\begin{equation}\label{Ribbonii}
  \Ribbon\VertexIII{s}{a}{t}{b}=
  \q^{\half((a + b - 2 t + 1)^2 - 1)}
  \sum_{i=0}^{s} \q^{-i a}
  \xi^i
  \Abin{t + i}{i}
  \prod_{j=0}^{i-1}\Aint{t + j - b}\,
  \VertexIII{s - i}{a}{t + i}{b}
\end{equation}
(we recall that $\xi=1-\q^2$).

\subsection{Algebra~\eqref{the-algebra} from~\eqref{the-loop}}
With the above ribbon structure, we now calculate
diagram~\eqref{YD-loop} in some cases.  To maintain association with
the diagram, we write $\chi_{\RepZ}(\RepY)$ as $\RepY\actsright\RepZ$
(the reasons for choosing the right action are purely
notational$/$graphical).  The calculations in what follows are based on
a formula for the double braiding: for two one-vertex modules, the
last diagram in~\eqref{B+B2} evaluates as
\begin{multline}\label{B2-VV}
   \mathsf{B}^2\left(V^{a}_s\tensor V^{b}_t\right)
   =\sum_{n=0}^{s + t}
   \sum_{i=n}^{s + t}\sum_{j=0}^{\min(i, t)}
  \q^{a b + 2 j (j - 1) + (i - n - 1) (i - n)  - 2 b j + a (n - 2 i - t)}
  \\[-6pt]
  \times{}\xi^{i - j} \Abin{i}{j}
  \Abin{s + t - j}{s} \Abin{s + t - n}{i - n}
  \prod_{\ell=0}^{i - j - 1}\!\Aint{\ell + j - b}
  \;\VertexIII{s + t - n}{a}{n}{b}.
\end{multline}

\subsubsection{}If $\RepY$ is irreducible, $\RepY\cong\RepX(r)_{\nu}$,
then $\RepX(r)_{\nu}\actsright\RepZ$
can only amount to multiplication by a number; indeed, we find that
\begin{equation*}
  \text{for all $x\in\RepXi{a}$
    with  $(a)_p\neq p-1$,}\quad
  x\actsright\RepXi{b} = \eigen{a,b}x,
\end{equation*}
where
\begin{equation*}
  \eigen{a,b}  =
  \frac{\q^{(a + 1)(b + 1)} - \q^{-(a + 1)(b + 1)}}{\q^{a + 1} -
    \q^{-a - 1}}.
\end{equation*}

It is instructive to reexpress this eigenvalue by indicating the
representation labels rather than the relevant coinvariants: for $a =
r' - 1 - \nu' p$ and $b = r - 1 - \nu p$,\pagebreak[3] we find that
$\RepX(r')_{\nu'}\actsright\RepX(r)_{\nu}$ amounts to multiplication
by
\begin{align*}
  \lambda(r',\nu';r,\nu) &= (-1)^{\nu'(r+1) + \nu r' + p \nu \nu'}\,
  \ffrac{\q^{r' r} - \q^{-r' r}}{\q^{r'} - \q^{-r'}}\\
  &=(-1)^{\nu'(r+1) + \nu r' + p \nu \nu'} \sum_{i=1}^{r} \q^{r' (r +
    1 - 2 i)}.
  \\
  \intertext{The last form is also applicable in the case where
    $r'=p$, and $\RepS(p)_{\nu'}\actsright\RepX(r)_{\nu}$ amounts to
    multiplication by}
  \lambda(p,\nu';r,\nu) &= (-1)^{(\nu' + 1)(r-1-\nu p)} r.
\end{align*}

For $\RepY=\RepV[r]_{\nu}$ in~\eqref{X-X-nu}, it may be worth noting
that the identity $\lambda(r',\nu';r,\nu)=\lambda(p-r',\nu'+1;r,\nu)$,
$1\leq r'\leq p-1$, explicitly shows that the action is the same on
both subquotients.  

\subsubsection{}Next, the action
$\RepP[r']_{\nu'}\actsright\RepX(r)_{\nu}$ has a diagonal piece, given
again by multiplication by $\lambda(r',\nu';r,\nu)$, and a nondiagonal
piece, mapping the top subquotient in
\begin{equation*}
  \raisebox{-1.95\baselineskip}{$\RepP[r']_{\nu'}\ {}={}\ $}
  \xymatrix@=12pt{
    &\modX(p-r')_{\nu'+1}
    \ar@/_12pt/[dl]
    \ar[dr]
    &\\
    \modX(r')_{\nu'}\ar[dr] &&\modX(r')_{\nu'+2}\ar@/_12pt/[dl]
    \\
    &\modX(p-r')_{\nu'+1}&
  }
\end{equation*}
into the bottom subquotient.  
Specifically, 
in terms of the ``top'' and ``bottom'' elements defined in~\eqref{uab}
and~\eqref{vab}, we have
\begin{equation*}
  u^{a, b}_t
  (1)\actsright\RepX(r) =
  \eigen{r', \nu'; r, \nu} u^{a, b}_t
  (1) +
  \mu(r', \nu'; r, \nu)
  v^{a, b}_t
  (r + 1),
\end{equation*}
where
\begin{multline*}
  \mu(r',\nu';r,\nu) = (-1)^{1 + \nu' r + \nu r' + p \nu' \nu}
  \ffrac{\q - \q^{-1}}{(\q^{r'} - \q^{-r'})^3}
  \\
  {}\times\Bigl((\q^{r' r} - \q^{-r' r})(\q^{r'} + \q^{-r'})
    - r(\q^{r' r} + \q^{-r' r})(\q^{r'} - \q^{-r'})\Bigr).
\end{multline*}
Because ${}\actsright\RepX(r)$ commutes with the $\Nich\!_p$ action
and coaction, and because $\RepPii{0}{a}{t}{b}$ is generated by the
$\Nich\!_p$ action and coaction from $u^{a, b}_t(1)$, the action of
$\RepXi{c}$ is thus defined on all of~$\RepPii{0}{a}{t}{b}$.


\subsubsection{}
Let $\fX(r)_{\nu}$ and $\fP(r)_{\nu}$ be the respective operations
$\actsright\RepX(r)_{\nu}$ and $\actsright\RepP(r)_{\nu}$.  We then
have relations~\eqref{the-algebra}, which are the fusion algebra
in~\cite{[FHST]}.

We see explicitly from the above formulas that
$\RepA_{\nu'}\actsright\RepX(r)_{\nu}$ depends on both $\nu'$ and
$\nu$ only modulo~2.

\section{Conclusion}
The construction of multivertex Yetter--Drinfeld $\Nich(X)$-modules
has a nice combinatorial flavor: elements of the braided space $X$
populate line intervals separated by ``vertex operators''---elements
of another braided space~$Y$, as
  $\xymatrix@R=6pt@C=70pt{
  \ar@*{[|(1.6)]}@{-}|(.15){\cross}
  |(.25){\punct}|(.4){\cross}|(.5){\cross}|(.65){\punct} |(.8){\cross}
  |(.90){\punct}[0,2]&&}\!\!\!\!$.
This construction and the $\Nich(X)$ action on such objects are
``universal'' in that they are formulated at the level of the braid
group algebra and work for any braiding.  However, even for diagonal
braiding, extracting information such as fusion from Nichols algebras
by direct calculation is problematic, except for rank $1$ (and
maybe~$2$).  Much greater promise is held by the program of finding
the modular group representation and then extracting the fusion from a
generalized Verlinde formula like the one in~\cite{[GT]}.
Importantly, those Nichols algebras that \textit{are} related to CFT{}
(and some certainly are, cf.~\cite{[2-boson]}) presumably carry an
$SL(2,\oZ)$ representation on the center of their Yetter--Drinfeld
category.

Going beyond Nichols algebras $\Nich(X)$ may also be interesting, and
is meaningful from the CFT{} standpoint: adding the divided powers
such as $F(p)$ in our $\Nich\!_p$ case, which are not in $\Nich(X)$
but do act on $\Nich\!_p$-modules, would yield a braided (and, in a
sense, ``one-sided'') analogue of the infinite-dimensional quantum
group that is Kazhdan--Lusztig-dual to logarithmic CFT{} models viewed
as Virasoro-symmetric theories~\cite{[BFGT],[BGT]}.

\subsubsection*{Acknowledgments} The content of Sec.~\ref{sec:Nich}
and Secs.~\bref{sec:screenings}--\bref{sec:dressedVO} is the joint
work with I.~Tip\-unin~\cite{[STbr]}.  I am grateful to
N.~Andruskiewitsch, T.~Creutzig, J.~Fjelstad, J.~Fuchs,
A.~Gainutdinov, I.~Heckenberger, D.~Ridout, I.~Runkel, C.~Schweigert,
I.~Tipunin, A.~Virelizier, and S.~Wood for the very useful
discussions.  A.~Virelizier also brought paper~\cite{[Brug]} to my
attention.  This paper was supported in part by the RFBR grant
10-01-00408 and the RFBR--CNRS grant 09-01-93105.

\mbox{}

\appendix
\section{Yetter--Drinfeld modules}\label{app:YD-axiom}
In the category of left--left module comodules over a braided Hopf
algebra $\Nich$, a Yetter--Drinfeld (also called ``crossed'')
module~\cite{[Sch-H-YD],[Besp-TMF],[Besp-next]} is a left module under
an action \ $\begin{tangle} \vstr{80}\lu
  \object{\raisebox{5.9pt}{\kern-4pt\tiny$\blacktriangleright$}}
\end{tangle}\ :\Nich\tensor\cY\to\cY$ and left comodule under a coaction
\ $\begin{tangle}\vstr{80}\ld\end{tangle}\ :\cY\to\Nich\tensor\cY$ such
that the axiom
\begin{equation}\label{yd-axiom}
  \begin{tangles}{l}
    \step\fobject{\Nich}\step[2]\fobject{\cY}\\
    \vstr{90}\cd\step\id\\
    \vstr{50}\id\step[2]\hx\\
    \vstr{90}\lu[2] \object{\raisebox{7pt}{\kern-4pt\tiny$\blacktriangleright$}}
    \step\hd\\
    \vstr{90}\ld[2]\step\ddh\\
    \vstr{50}\id\step[2]\hx\\
    \vstr{90}\cu\step\id
  \end{tangles}\ \ = \ \
  \begin{tangles}{l}
    \step\fobject{\Nich}\step[3]\fobject{\cY}\\
    \cd\step\ld\\
    \id\step[2]\hx\step\id\\
    \cu\step\lu\object{\raisebox{8.5pt}{\kern-4pt\tiny$\blacktriangleright$}}
  \end{tangles}
\end{equation}
holds.  The category $\BByd$ of Yetter--Drinfeld $\Nich$-modules is
monoidal and braided.  The action and coaction on a tensor product of
Yetter--Drinfeld modules are diagonal, respectively given by
\begin{equation}\label{act-coact-YD}
  \begin{tangles}{l}
    \hcd\step\id\step\id\\
    \id\step\hx\step\id\\
    \lu\object{\raisebox{8.5pt}{\kern-4pt\tiny$\blacktriangleright$}}
    \step
    \lu\object{\raisebox{8.5pt}{\kern-4pt\tiny$\blacktriangleright$}}
  \end{tangles}
  \qquad\text{and}\qquad
  \begin{tangles}{l}
      \ld\step\ld\\
      \id\step\hx\step\id\\
      \hcu\step\id\step\id
  \end{tangles}
\end{equation}
For two Yetter--Drinfeld modules, their braiding and its inverse and
square are given by
\begin{equation}\label{B+B2}
  \begin{tangles}{l}
    \vstr{67}\step[.5]\id\step\id\\
    \obox{2}{\mathsf{B}}\\    
    \vstr{67}\step[.5]\id\step\id
  \end{tangles}
  \ \ = \ \ 
  \begin{tangles}{l}
    \ld\step[1]\id\\
    \vstr{50}\id\step\hx\\
    \lu\object{\raisebox{8pt}{\kern-4pt\tiny$\blacktriangleright$}}\step[1]\id
  \end{tangles}
  \ ,
  \quad
  \begin{tangles}{l}
    \step[.5]\id\step\id\\
    \obox{2}{\mathsf{B}^{-1}}\\    
    \step[.5]\id\step\id
  \end{tangles}
  \ \ = \ \ 
  \begin{tangles}{l}
    \vstr{50}\step[1]\hx\\
    \ld\step[1]\id\\
    \vstr{67}\hxx\step[1]\id\\
    \id\step\O{\medA^{\scriptscriptstyle-1}}\step[1]\id\\
    \id\step\lu\object{\raisebox{8pt}{\kern-4pt\tiny$\blacktriangleright$}}
  \end{tangles}\ ,
  \quad\text{and}\quad
  \begin{tangles}{l}
    \step[.5]\id\step\id\\
    \obox{2}{\mathsf{B}^{2}}\\    
    \step[.5]\id\step\id
  \end{tangles}
  \ \ = \ \ 
  \begin{tangles}{l}
    \ld\step[2]\id\\
    \id\step\x\\
    \lu\object{\raisebox{8pt}{\kern-4pt\tiny$\blacktriangleright$}}\step[2]\id\\
    \ld\step[2]\id\\
    \id\step\x\\
    \lu\object{\raisebox{8pt}{\kern-4pt\tiny$\blacktriangleright$}}\step[2]\id
  \end{tangles}
  \ \ = \ \ 
  \begin{tangles}{l}
    \step[.5]\ld[2]\step[2]\id\\
    \hcd\step[1.5]\x\\
    \id\step\id\step[1]\hld\step[1]\ld\\
    \id\step\hx\step[.5]\id\step[1]\O{\medA}\step[1]\id\\[-.8pt]
    \hcu\step[1]\hlu\object{\raisebox{13.1pt}{\kern-4pt\tiny$\blacktriangleright$}}\step\lu\object{\raisebox{8pt}{\kern-4pt\tiny$\blacktriangleright$}}\\
    \step[.5]\id\step[2]\x\\
    \step[.5]\lu[2]\object{\raisebox{8pt}{\kern-4pt\tiny$\blacktriangleright$}}\step[2]\id
  \end{tangles}
\end{equation}

\section{Construction of Yetter--Drinfeld $\Nich\!_p$
  modules}\label{app:modules}

\subsection{One-vertex modules}\label{app:1-vertex}
One-vertex Yetter--Drinfeld $\Nich\!_p$-modules~\cite{[STbr]} are
spanned by $\VertexI{s}{a}$ (see~\eqref{Vsa}) for a fixed $a\in\oZ$
with $s$ ranging over a subset of $[0,\dots,p-1]$, under the
$\Nich\!_p$ action and coaction given in~\bref{sec:dressedVO}.

\subsubsection{Simple modules $\RepXi{a}$}\label{sec:Xi}
From each left coinvariant $\VertexI{0}{a}$, the action of $\Nich\!_p$
generates a simple module comodule of dimension $(a)_p+1$:
\begin{equation*}
  \RepXi{a}=\text{Span}(\VertexI{s}{a}\mid 0\leq s\leq (a)_p)
\end{equation*}
(simply because $F\adjoint\VertexI{(a)_p}{a}=0$ in accordance with the
above formulas).  The module comodule structure (in particular, the
matrix of $F(r)\adjoint$ in the basis of $\VertexI{s}{a}$,
Eq.~\eqref{Ad1-F(r)}) depends on $a$ only modulo~$p$, and hence there
are just $p$ nonisomorphic simple one-vertex module comodules, for
which we choose the notation $\RepX(r)$ indicating the dimension
$1\leq r\leq p$; then there are the $\Nich\!_p$ module comodule
isomorphisms
\begin{gather*}
  \RepXi{a}\cong\RepX(r)\quad
  \text{whenever} \quad (a)_p = r - 1.
\end{gather*}

\subsubsection{}\label{sec:S(p)}
As noted above, we sometimes use a special notation
$\RepS(p)=\RepX(p)$.

\subsubsection{}\label{sec:Xext}
For each $1\leq r\leq p-1$, $\RepXi{a}$ extends to a reducible module
comodule $\RepV[r]$ with $\RepXi{a}\cong\RepX(r)$ as a sub(co)module
and with the quotient isomorphic to $\RepX(p-r)$,\pagebreak[3] as
shown in~\eqref{V-def}.
In terms of basis, this is
\begin{equation}\label{X-X}
  \xymatrix@=15pt{
    &&&v^{a}(r+1)\ar@/_12pt/[dl]\ar^F[r]&\dots\ar^F[r]&v^{a}(p)\\
    v^{a}(1)\ar^F[r]&\dots\ar^F[r]&
    v^{a}(r)}
\end{equation}
where
\begin{equation*}
  v^a(i)= F^{i - 1}\adjoint \VertexI{0}{a},\quad i \leq r,
\end{equation*}
is a basis in $\RepXi{a}(r)$, with the last, $r$th element given by
$v^a(r)=C(r)\VertexI{r - 1}{a}$ with a nonzero $C(r)$,
and hence the upper floor starts with the element
$v^a(r+1)=C(r)\VertexI{r}{a}$.  The downward arrow in~\eqref{X-X} can
be understood to mean $X\mapsto x_1$ whenever $\delta X=\sum_s
F(s)\tensor x_s$; this convention is a reasonable alternative to
representing the same diagram as
\begin{equation*}
  \xymatrix@=15pt{
    &&&v^{a}(r+1)\ar^F[r]
    \ar@{{.}{--}{>}}|{F(r)}@/_18pt/[dlll]
    \ar@{{.}{--}{>}}@/_12pt/[dll]
    \ar@{{.}{--}{>}}|{F(1)}@/_8pt/[dl]
    &\dots\ar^F[r]&v^{a}(p)\\
    v^{a}(1)\ar^F[r]&\dots\ar^F[r]&
    v^{a}(r)}  
\end{equation*}
to express the idea that $\delta v^a(r+1)\in \Nich\!_p\tensor
\text{Span}(v^{a}(j)\mid 1\leq j\leq r + 1)$.

The general form of the adjoint action on the one-vertex space is
\begin{align}
  \label{Ad1-F(r)}
  F(r)\adjoint
  \VertexI{s}{a} &=
  \Abin{r + s}{r}\,\xi^r
  \prod_{i=s}^{s + r - 1}
  \Aint{i - a} \VertexI{s + r}{a},
\end{align}


\subsubsection{}
We verify that~\eqref{V1-decomp} holds by counting the total dimension
of the modules just constructed:
\begin{equation*}
  \dim\RepS(p)+\sum_{r=1}^{p-1}\dim\RepV[r]=p + (p - 1) p = p^2.
\end{equation*}

With the braiding~\eqref{Psi-F-V}, each of the above module comodules
satisfies the Yetter--Drinfeld axiom.

\subsection{Two-vertex modules}
A two-vertex Yetter--Drinfeld module is a linear span of some
$\VertexIII{s}{a}{t}{b}$, $0\leq s,t \leq p-1$, for fixed integers $a$
and $b$ (see~\eqref{Vsatb}).\pagebreak[3] The left adjoint action of
$\Nich\!_p$ on these is given by
\begin{align}
  \label{Ad2-F(r)}
  F(r)\adjoint
  \VertexIII{s}{a}{t}{b} &= \sum_{u=0}^{r}
  \coeff{a}{b}{s}{t}{r}{u}
  \VertexIII{s + r - u}{a}{t + u}{b},
\end{align}
where
\begin{equation}\label{c-st}
  \coeff{a}{b}{s}{t}{r}{u}=
  \xi^r
  \q^{u (2 s - a)}
  \abin{s + r - u}{r - u}\,
  \abin{t + u}{u}\,
  \prod_{i=u}^{r - 1}\Aint{s + i  + 2 t - a - b}\,
  \prod_{j=0}^{u - 1}\Aint{t + j - b}.
\end{equation}

The dependence on $b$ in~\eqref{c-st} is modulo~$p$, and on $a$,
modulo $2p$.  However, $\coeff{a+p}{b}{s}{t}{r}{u}=(-1)^u
\coeff{a}{b}{s}{t}{r}{u}$ and the matrix of $F(r)\adjoint$ in the
basis $(\VertexIII{s}{a+p}{t}{b})_{0\leq s,t\leq p-1}$ is the same as
in the basis $((-1)^t\VertexIII{s}{a}{t}{b})_{0\leq s,t\leq p-1}$;
moreover, the coaction is unaffected by this extra sign.  Hence, the
\textit{module comodule structure depends on both $a$ and $b$
  modulo~$p$}.

We arrive at decomposition~\eqref{V2-decomp} by first listing all
module comodules generated from \textit{left coinvariants}
\begin{equation*}
   \VertexIII{0}{a}{t}{b}=
   \xymatrix@1@C=70pt{
    \ar@*{[|(1.6)]}@{-}|(.20){\punct}^(.20){{}^{\scriptstyle
        a}}|(.3){\cross}|(.50){\cross}|(.65){\cross}|(.85){\punct}^(.85){{}^{\scriptstyle
        b}}[0,2]&&} \qquad\text{($t$ crosses),}
\end{equation*}

\noindent
and then studying their extensions.

In accordance with~\eqref{Ad2-F(r)}, the algebra acts on left
coinvariants as
\begin{equation*}
  F(r)\adjoint
  \VertexIII{0}{a}{t}{b} = \sum_{s=0}^{r}
  c^{a,b}_t(r,s)
  \VertexIII{r - s}{a}{t + s}{b},
\end{equation*}
with the coefficients
\begin{equation}\label{c-t}
  c^{a,b}_t(r,s)=
  c^{a,b}_{0,t}(r,s)=
  \xi^r
  \q^{-s a} \abin{t + s}{s}
  \prod_{i=s}^{r - 1}\Aint{i + 2 t - a - b}\,
  \prod_{j=0}^{s - 1}\Aint{t + j - b}. 
\end{equation}
In practical terms, the cases in~\bref{sec:all-cases} can be
conveniently studies as follows.
\begin{enumerate}\addtolength{\itemsep}{4pt}
  \renewcommand{\labelenumi}{(\textbf{\arabic{enumi}})}
\renewcommand{\theenumi}{(\textbf{\arabic{enumi}})}

\item\label{case:Steinberg} $F(r)\adjoint \VertexIII{0}{a}{t}{b}$ is
  nonvanishing and not a coinvariant for all $r$, $1\leq r\leq p-1$.

\item\label{case:short-vanishing2}$\VertexIII{0}{a}{t}{b}$ is not in
  the image of $F$ and $F(r)\adjoint\VertexIII{0}{a}{t}{b}$ vanishes
  for some $r\leq p-1$, i.e.,
  \begin{equation*}
    c^{a,b}_t(r,s)=0, \quad 1\leq s\leq r.
  \end{equation*}

\item\label{case:coinvariant} $F(r)\adjoint \VertexIII{0}{a}{t}{b}$ is
  a coinvariant, i.e., $F(r)\adjoint
  \VertexIII{0}{a}{t}{b}=\mathrm{const}\;\VertexIII{0}{a}{t+r}{b}$,
  for some $r\leq p-1$, which is equivalent to
  \begin{equation*}
    \begin{cases}
      c^{a,b}_t(r,s)=0,&0\leq s\leq r-1,\\
      c^{a,b}_t(r,r)\neq 0.
    \end{cases}
  \end{equation*}

\end{enumerate}


Let $\beta=(a+b-2t)_p + 1$.  For $0\leq b\leq p-1$, this is equivalent
to $b=(2t + \beta - 1 - a)_p$.  In fact, every triple $(a,b,t)$,
$0\leq a,b,t\leq p-1$, can be uniquely represented as
\begin{equation}\label{abt-param}
  (a,b,t)=(a,\; (2 t + \beta - 1 - a)_p,\; t),
  \quad
  1\leq \beta\leq p.
\end{equation}
In this parameterization, coefficients~\eqref{c-t} become
\begin{gather*}
  c^{a,b}_t(r,s) =
  \xi^r
  \q^{-s a} \abin{t + s}{s}
  \prod_{i=s+1}^{r}\Aint{i - \beta}\,
  \prod_{j=1}^{s}\Aint{j - \beta + a - t}.
\end{gather*}
and the analysis of the above cases becomes relatively
straightforward.  The results are as
follows.\enlargethispage{\baselineskip}

\subsubsection{Irreducible dimension-$p$ modules $\RepSii{0}{a}{t}{b}(p)$
  (case~\ref{case:Steinberg})}\label{sec:S}
A simple module comodule of dimension $p$, isomorphic to $\RepS(p)$
in~\bref{sec:S(p)}, is generated under the $\Nich\!_p$ action from a
coinvariant $\VertexIII{0}{a}t{b}$ if and only if
\begin{equation}\label{cond:S}
  (a + b - 2 t)_p + 1=p.\pagebreak[3]
\end{equation}\pagebreak[3]
When this condition is satisfied, we write $\RepSii{0}{a}t{b}$, or even
$\RepSii{0}{a}t{b}(p)$, for this module comodule isomorphic to
$\RepS(p)$.\footnote{Condition~\eqref{cond:S} is actually worked out
  as follows: For odd $p$, it holds if and only if either $t \equiv
  \Mod{\half(a + b + 1 + p)}{p}$ with $a + b$ even, or $t = \half(a +
  b + 1)$ with $a + b$ odd.  For even $p$, it holds if and only if
  either $t \equiv\Mod{\half(p + 1 + a + b)}{p}$ or $t = \half(a + b +
  1)$ (which selects only odd $a+b$).}

\subsubsection{Reducible dimension-$p$ modules
  $\RepVii{0}{a}{t}{b}[r]$
  (case~\ref{case:short-vanishing2})}\label{sec:X}
A simple module comodule isomorphic to $\RepX(r)$ for some $1\leq
r\leq p-1$ is generated under the action of $\Nich\!_p$ from a
coinvariant $\VertexIII{0}{a}t{b}$ that is not itself in the image of
$F$ if and only if $r=(a + b - 2 t)_p + 1$ and either of the two
conditions holds:\footnote{The logic of the presentation is that we
  assume that $0\leq a,b\leq p-1$, and hence $(a)_p=a$; but we do not
  omit the operator of taking the residue modulo~$p$ because we refer
  to formulas given here also in the case where $a\in\oZ$.}
\begin{align}
  \label{cond:Xi}
  &t\leq (a)_p\quad\text{and} \quad (a)_p - r + 1\leq t\leq p - 1 - r,\\  
  \label{cond:Xii}
  &t\geq (a)_p + 1\quad\text{and} \quad
  p - r\leq t\leq p - r + (a)_p.
\end{align}
In this case, we write $\RepXii{0}{a}t{b}$ or $\RepXii{0}{a}t{b}(r)$
for the corresponding module comodule:
\begin{equation}\label{2X-mod}
  \RepXii{0}{a}t{b}\cong\RepX(r)\quad\text{whenever \
    $r=(a + b - 2 t)_p + 1$ and~\eqref{cond:Xi}
    or~\eqref{cond:Xii} holds.}
\end{equation}

Every $\RepXii{0}{a}t{b}$ is further extended as in~\eqref{V-def},
which in terms of basis is now realized as
\begin{equation*}
  \xymatrix@=15pt@C=12pt{
    &&&*{\sum\limits_{s=0}^{r-1}
    \Afac{r\!-\!1}c^{a,b}_t(r\!-\!1,s)\!\VertexIII{r - s}{a}{t + s}{b}}
    \ar@/_20pt/[dl]\ar^(.8)F[r]
    &\dots
    \\
    \VertexIII{0}{a}t{b}\!\ar^(.65)F[r]&\dots\ar^(.17)F[r]&
    *{\sum\limits_{s=0}^{r-1}
    \Afac{r\!-\!1}c^{a,b}_t(r\!-\!1,s) \VertexIII{r -1 - s}{a}{t + s}{b}}
    }
\end{equation*}
with the south-west arrow meaning the same as in~\bref{sec:Xext}; the
quotient is isomorphic to $\RepX(p-r)$.  The notation\pagebreak[3]
$\RepVii{0}{a}{t}{b}[r]$ for this dimension-$p$ module comodule
explicitly indicates the relevant left coinvariant and the dimension
of the sub(co)module; the module comodule structure depends only on
$r$: $\RepVii{0}{a}{t}{b}[r]\cong\RepV[r]$.

\subsubsection{Three-floor modules $\RepPii{0}{a}{t}{b}[r]$
  (case~\ref{case:coinvariant})}\label{sec:LB}
We next assume that none of the above conditions \eqref{cond:S},
\eqref{cond:Xi}, \eqref{cond:Xii} is satisfied.  An exemplary
exercise shows that the negation of
$\eqref{cond:S}\lor\eqref{cond:Xi}\lor\eqref{cond:Xii}$ is the ``or''
of the four conditions
\begin{align}
  \label{P-cond-3}
  &t \geq p - r + (a)_p + 1,\\
  \label{P-cond-4}
  &p - r \leq t \leq (a)_p,\\
  \label{P-cond-1}
  &t \leq (a)_p - r,\\
  \label{P-cond-2}
  &(a)_p + 1 \leq t \leq p - r - 1,
\end{align}
where again $r=(a + b - 2 t)_p + 1$, $1\leq r\leq p-1$.  The module
generated from the coinvariant $\VertexIII{0}{a}t{b}$ is then a
sub(co)module in an indecomposable module comodule with the structure
of subquotients
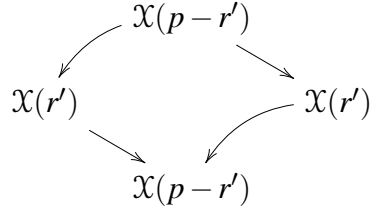
\begin{equation*}
  \xymatrix@=12pt{
    &&\modX(p-r')
    \ar@/_12pt/[dl]
    \ar[dr]
    &\\
    &\modX(r')\ar[dr]
    &
    &\modX(r')\ar@/_12pt/[dl]
    \\
    &&\modX(p-r')&
  }
\end{equation*}
where $r'$ is either $r$ or $p-r$, as we now
describe.

\begin{enumerate}\addtolength{\itemsep}{4pt}%
  \renewcommand{\labelenumi}{\roman{enumi}.}%
\renewcommand{\theenumi}{\roman{enumi}}%
\item\label{case:bot}If $t+r\geq p$ (which means that either
  \eqref{P-cond-3} or \eqref{P-cond-4} holds), then the submodule
  generated from $\VertexIII{0}{a}t{b}$ is isomorphic to $\RepX(r)$.
  We let it be denoted by $\RepBii{0}{a}t{b}(r)$. ($\RepB$~is for
  ``bottom,'' and $\RepL$ is for ``left.'')

\item\label{case:left}If $t+r\leq p-1$ (which means that either
  \eqref{P-cond-1} or \eqref{P-cond-2} holds), then the submodule
  generated from $\VertexIII{0}{a}t{b}$, denoted by
  $\RepLii{0}{a}t{b}(r)$, is a $p$-dimensional reducible module
  comodule with $\RepBii{0}{a}{t+r}{b}(p-r)\cong\RepX(p-r)$ as a
  sub(co)module and with the quotient isomorphic to $\RepX(r)$:
  \begin{equation*}
    \raisebox{-\baselineskip}{$\RepLii{0}{a}t{b}[r]={}$\ \ \ }
    \xymatrix@=12pt{\RepX(r)\ar[dr]\\
      &
      \RepX(p-r).
    }
  \end{equation*}
\end{enumerate}
In terms of basis, this diagram is
\begin{equation*}
  \xymatrix@C=15pt@R=12pt{
    \VertexIII{0}{a}{t}{b}\ar^F[r]&\dots\ar^(.15)F[r]&
    \sum\limits_{s=0}^{r-1} \Afac{r-1}c^{a,b}_t(r-1,s)
    \VertexIII{r -1 - s}{a}{t + s}{b}\ar[dr]^(.6)F\\
    &&&\Afac{r}\VertexIII{0}{a}{t+r}{b}\ar^(.55)F[r]&\dots
  }
\end{equation*}

The set of all diagrams of this type actually describes both
cases~\ref{case:bot} and \ref{case:left}: according to whether a given
coinvariant $\VertexIII{0}{a}{u}{b}$ is or is not in the image of $F$,
it occurs either in the bottom line (case \ref{case:bot}) or in the
upper line (case \ref{case:left}) of the last diagram.

Every such diagram is extended further, again simply because of the
``cofree'' nature of the coaction:
\begin{gather}\label{Xtended}
  {}\kern-20pt\xymatrix@C=15pt@R=12pt{
    &&& *{\ T^{a,b}_{t+r}(r)} \ar@/_12pt/[ld]
    \\
    \VertexIII{0}{a}{t}{b}\ar^F[r]&\dots\ar^(.15)F[r]&
    *{\sum\limits_{s=0}^{r-1}\Afac{r-1} c^{a,b}_t(r-1,s)
      \VertexIII{r - s - 1}{a}{t + s}{b}}
    \ar[dr]^(.6)F\\
    &&&*{\Afac{r}c^{a,b}_t(r,r)\VertexIII{0}{a}{t+r}{b}\ }
    \ar^(.75)F[r]&\dots
    }\kern-20pt
\end{gather}
where, evidently,
\begin{equation*}
  T^{a,b}_{t+r}(r)=
  \sum_{s=0}^{r - 1}
  \Afac{r-1}c^{a, b}_t(r - 1, s) \VertexIII{r - s}{a}{t + s}{b}.
\end{equation*}
Setting
\begin{align}
  \label{uab}
  u^{a, b}_t(i)&=F^{i-1}\adjoint T^{a,b}_{t+r}(r),\\
  \label{vab}
  v^{a, b}_t(i)&=F^{i-1}\adjoint\VertexIII{0}{a}{t}{b},
\end{align}
we have the full picture extending~\eqref{Xtended} as (omitting the
${}^{a,b}_{t}$ labels for brevity)
\begin{equation*}
  \xymatrix@C=12pt@R=12pt{
    &&& u(1) \ar@/_12pt/[dl]\ar^(.55)F[r]&\dots\ar^(.45)F[r]&
    u(p-r)\ar^F[dr]
    \\
    v(1)\ar^F[r]&\dots\ar^(.45)F[r]&
    v(r)\ar[dr]^(.55)F&&&&u(p-r+1)\ar@/_12pt/[ld]
    \ar^(.65)F[r]&\dots\ar^(.45)F[r]&u(p)\\
    &&&v(r+1)
    \ar^(.55)F[r]&\dots\ar^(.55)F[r]&v(p)
    }
\end{equation*}
Here,\footnote{``The closure of the rhombus'' in the above diagram is
  a good illustration of the use of the Yetter--Drinfeld axiom, which
  is also used in several other derivations without special notice.
  The ``relative factor'' $\q^{-2r}$ in the next two formulas, in
  particular, is an immediate consequence of the Yetter--Drinfeld
  condition.}
\begin{equation*}
  \delta u(1)
  =
  1\tensor u(1) +
  \ffrac{1}{\Afac{r-1}}\,F(1)\tensor v(r)+\dots
\end{equation*}
and, similarly,
\begin{equation*}
  \delta u(p+r-1) = 1\tensor u(p+r-1) +
  \ffrac{\q^{-2 r}}{\Afac{r - 1}}\,F(1)\tensor v(p) + \dots.
\end{equation*}

To label such modules by the leftmost coinvariant
$\VertexIII{0}{a}{t}{b}$ (even though the entire module is not
\textit{generated from} this element), we write $\RepPii{0}{a}{t}{b}$
to indicate both the module type and the characteristic coinvariant.
An even more redundant notation is
$\RepPii{0}{a}{t}{b}[r]$, indicating the length $r$ of the left wing
(which of course is $r=(a+b-2t)_p+1$).  The module comodule structure
depends only on~$r$:
\begin{equation}\label{P-iso}
  \RepPii{0}{a}{t}{b}[r] \cong \RepP[r].
\end{equation}
To summarize,
given a coinvariant $\VertexIII{0}{a}{t}{b}$, \eqref{P-iso} holds if
and only if (for $r=(a + b - 2 t)_p + 1$)
\begin{equation}\label{ProjectiveQ}
  \begin{aligned}
    &1\leq r\leq p-1 \text{ \ and}\\
    &\text{$t \leq (a)_p - r$ \ or \
      $(a)_p + 1 \leq t \leq p - r - 1$.}
  \end{aligned}  
\end{equation}

\subsubsection{Completeness}We verify~\eqref{V2-decomp} by counting
the total dimension of the modules constructed.  This gives $p^4$, the
dimension of $\mathbb{V}_p(2)$, as follows.\pagebreak[3] There are
$p^2$ modules $\RepS(p)$ constructed in~\bref{sec:S}, $2 r (p - r)$
modules $\RepX(r)$ in~\bref{sec:X} for each $1\leq r\leq p-1$, making
the total of $\frac{1}{3} p (p^2 - 1)$, and, finally, $(p - r)^2$
modules $\RepL(r)$ in~\bref{sec:LB} for each $1\leq r\leq p-1$, making
the total of $\frac{1}{6} p (p - 1) (2 p - 1)$.  Each $\RepS(p)$ is
$p$-dimensional, each $\RepX(r)$ \textit{extends to} a $p$-dimensional
module, and each $\RepLii{0}{a}t{b}[r]$ extends to a $2p$-dimensional
module.
The total dimension is
\begin{equation*}
  p^2\cdot p + \ffrac{1}{3}\, p (p^2 - 1) \cdot p +
  \ffrac{1}{6}\,p (p - 1)(2 p - 1)\cdot 2 p = p^4.
\end{equation*}


\subsubsection{Example} Decomposition~\eqref{V2-decomp} is illustrated
in Fig.~\ref{fig:5-modules} for $p=5$.
\begin{figure}[tbp]
  \centering\footnotesize
  \begin{equation*}
    \begin{matrix}
      \begin{pmatrix}
        \RepXii{0}{0}{0}{0}\left(1\right)_0 \\
        \RepXii{0}{0}{1}{0}\left(4\right)_1 \\
        \RepLii{0}{0}{2}{0}\left[2\right]_1 \\
        \RepSii{0}{0}{3}{0}\left(5\right)_2 \\
        \RepBii{0}{0}{4}{0}\left(3\right)_2
      \end{pmatrix}
      &\begin{pmatrix}
        \RepXii{0}{0}{0}{1}\left(2\right)_0 \\
        \RepSii{0}{0}{1}{1}\left(5\right)_1 \\
        \RepXii{0}{0}{2}{1}\left(3\right)_1 \\
        \RepLii{0}{0}{3}{1}\left[1\right]_1 \\
        \RepBii{0}{0}{4}{1}\left(4\right)_2
      \end{pmatrix}
      &\begin{pmatrix}
        \RepXii{0}{0}{0}{2}\left(3\right)_0 \\
        \RepLii{0}{0}{1}{2}\left[1\right]_0 \\
        \RepBii{0}{0}{2}{2}\left(4\right)_1 \\
        \RepXii{0}{0}{3}{2}\left(2\right)_1 \\
        \RepSii{0}{0}{4}{2}\left(5\right)_2
      \end{pmatrix}
      &\begin{pmatrix}
        \RepXii{0}{0}{0}{3}\left(4\right)_0 \\
        \RepLii{0}{0}{1}{3}\left[2\right]_0 \\
        \RepSii{0}{0}{2}{3}\left(5\right)_1 \\
        \RepBii{0}{0}{3}{3}\left(3\right)_1 \\
        \RepXii{0}{0}{4}{3}\left(1\right)_1
      \end{pmatrix}
      &\begin{pmatrix}
        \RepSii{0}{0}{0}{4}\left(5\right)_0 \\
        \RepLii{0}{0}{1}{4}\left[3\right]_0 \\
        \RepLii{0}{0}{2}{4}\left[1\right]_0 \\
        \RepBii{0}{0}{3}{4}\left(4\right)_1 \\
        \RepBii{0}{0}{4}{4}\left(2\right)_1
      \end{pmatrix}
      \\[8pt]
      \begin{pmatrix}
        \RepXii{0}{1}{0}{0}\left(2\right)_0 \\
        \RepSii{0}{1}{1}{0}\left(5\right)_1 \\
        \RepXii{0}{1}{2}{0}\left(3\right)_1 \\
        \RepLii{0}{1}{3}{0}\left[1\right]_1 \\
        \RepBii{0}{1}{4}{0}\left(4\right)_2
      \end{pmatrix}
      &\begin{pmatrix}
        \RepXii{0}{1}{0}{1}\left(3\right)_0 \\
        \RepXii{0}{1}{1}{1}\left(1\right)_0 \\
        \RepXii{0}{1}{2}{1}\left(4\right)_1 \\
        \RepXii{0}{1}{3}{1}\left(2\right)_1 \\
        \RepSii{0}{1}{4}{1}\left(5\right)_2
      \end{pmatrix}
      &\begin{pmatrix}
        \RepXii{0}{1}{0}{2}\left(4\right)_0 \\
        \RepXii{0}{1}{1}{2}\left(2\right)_0 \\
        \RepSii{0}{1}{2}{2}\left(5\right)_1 \\
        \RepXii{0}{1}{3}{2}\left(3\right)_1 \\
        \RepXii{0}{1}{4}{2}\left(1\right)_1
      \end{pmatrix}
      &\begin{pmatrix}
        \RepSii{0}{1}{0}{3}\left(5\right)_0 \\
        \RepXii{0}{1}{1}{3}\left(3\right)_0 \\
        \RepLii{0}{1}{2}{3}\left[1\right]_0 \\
        \RepBii{0}{1}{3}{3}\left(4\right)_1 \\
        \RepXii{0}{1}{4}{3}\left(2\right)_1
      \end{pmatrix}
      &\begin{pmatrix}
        \RepLii{0}{1}{0}{4}\left[1\right]_{-1} \\
        \RepBii{0}{1}{1}{4}\left(4\right)_0 \\
        \RepLii{0}{1}{2}{4}\left[2\right]_0 \\
        \RepSii{0}{1}{3}{4}\left(5\right)_1 \\
        \RepBii{0}{1}{4}{4}\left(3\right)_1
      \end{pmatrix}
      \\
      \hdotsfor{5}
      \\
      \begin{pmatrix}
        \RepSii{0}{4}{0}{0}\left(5\right)_0 \\
        \RepLii{0}{4}{1}{0}\left[3\right]_0 \\
        \RepLii{0}{4}{2}{0}\left[1\right]_0 \\
        \RepBii{0}{4}{3}{0}\left(4\right)_1 \\
        \RepBii{0}{4}{4}{0}\left(2\right)_1
      \end{pmatrix}
      &\begin{pmatrix}
        \RepLii{0}{4}{0}{1}\left[1\right]_{-1} \\
        \RepBii{0}{4}{1}{1}\left(4\right)_0 \\
        \RepLii{0}{4}{2}{1}\left[2\right]_0 \\
        \RepSii{0}{4}{3}{1}\left(5\right)_1 \\
        \RepBii{0}{4}{4}{1}\left(3\right)_1
      \end{pmatrix}
      &\begin{pmatrix}
        \RepLii{0}{4}{0}{2}\left[2\right]_{-1} \\
        \RepSii{0}{4}{1}{2}\left(5\right)_0 \\
        \RepBii{0}{4}{2}{2}\left(3\right)_0 \\
        \RepLii{0}{4}{3}{2}\left[1\right]_0 \\
        \RepBii{0}{4}{4}{2}\left(4\right)_1
      \end{pmatrix}
      &\begin{pmatrix}
        \RepLii{0}{4}{0}{3}\left[3\right]_{-1} \\
        \RepLii{0}{4}{1}{3}\left[1\right]_{-1} \\
        \RepBii{0}{4}{2}{3}\left(4\right)_0 \\
        \RepBii{0}{4}{3}{3}\left(2\right)_0 \\
        \RepSii{0}{4}{4}{3}\left(5\right)_1
      \end{pmatrix}
      &\begin{pmatrix}
        \RepLii{0}{4}{0}{4}\left[4\right]_{-1} \\
        \RepLii{0}{4}{1}{4}\left[2\right]_{-1} \\
        \RepSii{0}{4}{2}{4}\left(5\right)_0 \\
        \RepBii{0}{4}{3}{4}\left(3\right)_0 \\
        \RepBii{0}{4}{4}{4}\left(1\right)_0
      \end{pmatrix}
    \end{matrix}
  \end{equation*}
  \caption{\small For each $0\leq t\leq p-1$, $0\leq a\leq p-1$,
    $0\leq b\leq p-1$ (where $p=5$), the module comodule
    \textit{generated} from $\VertexIII{0}{a}{t}{b}(r)$ is indicated
    as $\modA^{a,\,b}_{0,\,t}(r)_{\nu}$, where $r$ is the dimension of
    the relevant subquotient, $\nu$ is the braiding index, and $\modA$
    indicates the module type.  Only $a=0,1,4$ are shown for
    compactness.  Whenever an $\RepLii{0}{a}{t}{b}[r]$ occurs in a
    column of height~5, the $\RepBii{0}{a}{t+r}{b}(p-r)$ module is
    present in the same column.  We do not replace negative braiding
    indices $-1$ with the ``canonical'' representative 3 in $\oZ_4$
    ``for continuity.''}
  \label{fig:5-modules}
\end{figure}
The figure lists all the modules generated from the
$\VertexIII{0}{a}{t}{b}$ with $a=0,1,4$, $b=0,1,2,3,4$, and
$t=0,1,2,3,4$ (two values of $a$ are omitted for compactness).  Each
$\RepBii{0}{a}{u}{b}(r)$ module is a Yetter--Drinfeld sumbodule in the
$\RepLii{0}{a}{u-r}{b}[p-r]$ module in the same column of $p=5$
modules.  The subscript additionally indicates the braiding sectors
(see~\bref{sec:braiding-sectors}).


\begin{thebibliography}{99}
\bibitem{[KLx]} D.\;Kazhdan and G.\;Lusztig, \textit{Tensor structures
    arising from affine Lie algebras,} I, J. Amer. Math. Soc. 6 (1993)
  905--947; II, J. Amer. Math. Soc. 6 (1993) 949--1011; III, J.\
  Amer.\ Math.\ Soc.\ 7 (1994) 335--381; IV, J.\ Amer.\ Math.\ Soc.\ 7
  (1994) 383--453.

\bibitem{[Fink]}M.\;Finkelberg, \textit{An equivalence of fusion
    categories}, Geometric and Functional Analysis (GAFA) 6 (1996)
  249--267.

\bibitem{[T]}V.G.\;Turaev, \textsl{Quantum Invariants of Knots and
    $3$-Manifolds}, Walter de~Gruyter, Berlin--New~York (1994).

\bibitem{[MS]}G.\;Moore and N.\;Seiberg, \textit{Lectures on RCFT}, in:
  \textsl{Physics, Geometry, and Topology} (Trieste spring school
  1989), p.~263; Plenum (1990).
 
\bibitem{[FRS]}J.\;Fuchs, I.\;Runkel, and C.\;Schweigert, \textit{TFT
    construction of RCFT correlators I: Partition functions},
  Nucl. Phys. B 646 (2002) 353 [hep-th$/$0204148]; \textit{TFT
    construction of RCFT correlators II: Unoriented surfaces},
  Nucl. Phys. B 678 (2004) 511 [hep-th$/$0306164]; \textit{TFT
    construction of RCFT correlators III: Simple currents},
  Nucl. Phys. B 694 (2004) 277 [hep-th$/$0403158]; \textit{TFT
    construction of RCFT correlators IV: Structure constants and
    correlation functions}, Nucl. Phys. B 715 (2005) 539
  [hep-th$/$0412290].

\bibitem{[FFRS]}J.\;Fr\"ohlich, J.\;Fuchs, I.\;Runkel, and
  C.\;Schweigert, \textit{Correspondences of ribbon categories},
  Adv. Math. 199 (2006) 192 [math.CT/0309465].

\bibitem{[FGST]}B.L.\;Feigin, A.M.\;Gainutdinov, A.M.\;Semikhatov, and
  I.Yu.\;Tipunin, \textit{Modular group representations and fusion in
    logarithmic conformal field theories and in the quantum group
    center}, Commun. Math. Phys. 265 (2006) 47--93
  [arXiv:\linebreak[0]hep-th/\linebreak[0]0504093].

\bibitem{[FGST2]}B.L.\;Feigin, A.M.\;Gainutdinov, A.M.\;Semikhatov, and
  I.Yu.\;Tipunin, \textit{Kazhdan--Lusztig correspondence for the
    representation category of the triplet $W$-algebra in logarithmic
    CFT}, Theor. Math. Phys. 148 (2006) 1210--1235
  [arXiv:\linebreak[0]math/0512621 [math.QA]].

\bibitem{[FGST3]}B.L.\;Feigin, A.M.\;Gainutdinov, A.M.\;Semikhatov, and
  I.Yu.\;Tipunin, \textit{Logarithmic extensions of minimal models:
    characters and modular transformations}, Nucl. Phys.  B757 (2006)
  303--343 [arXiv:\linebreak[0]hep-th/\linebreak[0]0606196].

\bibitem{[FGSTq]}B.L.\;Feigin, A.M.\;Gainutdinov, A.M.\;Semikhatov,
  and I.Yu.\;Tipunin, \textit{Kazhdan--Lusztig-dual quantum group for
    logarithmic extensions of Virasoro minimal models},
  J. Math. Phys. 48 (2007) 032303
  [arXiv:\linebreak[0]math/\linebreak[0]0606506 [math.QA]].

\bibitem{[NT]} K.\;Nagatomo and A.\;Tsuchiya, \textit{The triplet
    vertex operator algebra $W(p)$ and the restricted quantum group at
    root of unity}, arXiv:\linebreak[0]0902.4607 [math.QA].

\bibitem{[GT]} A.M.\;Gainutdinov and I.Yu.\;Tipunin, \textit{Radford,
    Drinfeld, and Cardy boundary states in $(1,p)$ logarithmic
    conformal field models}, J. Phys. A42 (2009) 315207
  [arXiv:\linebreak[0]0711.3430].

\bibitem{[BFGT]} P.V.\;Bushlanov, B.L.\;Feigin, A.M.\;Gainutdinov, and
  I.Yu.\;Tipunin, \textit{Lusztig limit of quantum $sl(2)$ at root of
    unity and fusion of $(1,p)$ Virasoro logarithmic minimal models},
  Nucl. Phys. B818 (2009) 179--195 [arXiv:\linebreak[0]0901.1602].

\bibitem{[BGT]} P.V.\;Bushlanov, A.M.\;Gainutdinov, and
  I.Yu.\;Tipunin, \textit{Kazhdan-Lusztig equivalence and fusion of
    Kac modules in Virasoro logarithmic models},
  arXiv:\linebreak[0]1102.0271.

\bibitem{[AA]}T.\;Abe and Y.\;Arike, \textit{Intertwining operators and
    fusion rules for vertex operator algebras arising from symplectic
    fermions}, arXiv:\linebreak[0]1108.1823.

\bibitem{[CR]} T.\;Creutzig and D.\;Ridout, \textit{Relating the
    archetypes of logarithmic conformal field theory}
  arXiv:\linebreak[0]1107.\linebreak[0]2135.

\bibitem{[HY]} Y.-Z.\;Huang and J.\;Yang, \textit{Logarithmic
    intertwining operators and associative algebras}
  arXiv:\linebreak[0]1104.\linebreak[0]4679.

\bibitem{[AN]}Y.\;Arike and K.\;Nagatomo, \textit{Some remarks on
    pseudo-trace functions for orbifold models associated with
    symplectic fermions}, arXiv:\linebreak[0]1104.0068.

\bibitem{[VJS]} R.\;Vasseur, J.L.\;Jacobsen, and H.\;Saleur,
  \textit{Indecomposability parameters in chiral logarithmic conformal
    field theory}, Nucl. Phys. B851 (2011) 314--345
  [arXiv:\linebreak[0]1103.3134].

\bibitem{[AM-2p]}D.\;Adamovic and A.\;Milas, \textit{On $W$-algebra
    extensions of $(2,p)$ minimal models: $p > 3$},
  arXiv:\linebreak[0]1101.\linebreak[0]0803.

\bibitem{[HLZ]}Y.-Z.\;Huang, J.\;Lepowsky, and L.\;Zhang,
  \textit{Logarithmic tensor category theory for generalized modules
    for a conformal vertex algebra, I: Introduction and strongly
    graded algebras and their generalized modules},
  arXiv:\linebreak[0]1012.4193; \textit{Logarithmic tensor category
    theory, II: Logarithmic formal calculus and properties of
    logarithmic intertwining operators}, arXiv:\linebreak[0]1012.4196;
  \textit{Logarithmic tensor category theory, III: Intertwining maps
    and tensor product bifunctors}, arXiv:\linebreak[0]1012.4197;
  \textit{Logarithmic tensor category theory, IV: Constructions of
    tensor product bifunctors and the compatibility conditions};
  arXiv:\linebreak[0]1012.4198; \textit{Logarithmic tensor category
    theory, V: Convergence condition for intertwining maps and the
    corresponding compatibility condition},
  arXiv:\linebreak[0]1012.4199; \textit{Logarithmic tensor category
    theory, VI: Expansion condition, associativity of logarithmic
    intertwining operators, and the associativity isomorphisms},
  arXiv:\linebreak[0]1012.4202.

\bibitem{[GRW]}M.R.\;Gaberdiel, I.\;Runkel, and S.\;Wood, \textit{A
    modular invariant bulk theory for the $c=0$ triplet model}
  J. Phys. A44 (2011) 015204 [arXiv:\linebreak[0]1008.0082].

\bibitem{[FSS-11]}J.\;Fuchs, C.\;Schweigert,, and C.\;Stigner,
  \textit{Modular invariant Frobenius algebras from ribbon Hopf
    algebra automorphisms}, arXiv:1106.0210.

\bibitem{[FSS-12]}J.\;Fuchs, C.\;Schweigert,, and C.\;Stigner,
  \textit{The Cardy--Cartan modular invariant}, arXiv:1201.4267.

\bibitem{[VGJS]} R.\;Vasseur, A.M.\;Gainutdinov, J.L.\;Jacobsen, and
  H.\;Saleur, \textit{The puzzle of bulk conformal field theories at
    central charge $c=0$}, arXiv:1110.1327.

\bibitem{[GV]} A.M.\;Gainutdinov and R.\;Vasseur, \textit{Lattice
    fusion rules and logarithmic operator product expansions},
  arXiv:1203.6289.

\bibitem{[RGW-12]}I.\;Runkel, M.R.\;Gaberdiel, and S.\;Wood,
  \textit{Logarithmic bulk and boundary conformal field theory and the
    full centre construction}, arXiv:1201.6273.

\bibitem{[STbr]} A.M.\;Semikhatov and I.Yu.\;Tipunin, \textit{The
    Nichols algebra of screenings}, arXiv:\linebreak[0]1101.5810.

\bibitem{[Nich]}W.\;D.\;Nichols, \textit{Bialgebras of type one},
  Commun. Algebra 6 (1978) 1521--1552.

\bibitem{[AG]}N.\;Andruskiewitsch and M.\;Gra\~na, \textit{Braided
    Hopf algebras over non abelian finite groups}, Bol. Acad. Nacional
  de Ciencias (Cordoba) 63 (1999) 45--78
  [arXiv:\linebreak[0]math$/$9802074 [math.QA]].

\bibitem{[AS-onthe]}N.\;Andruskiewitsch and H.-J.\;Schneider,
  \textit{On the classification of finite-dimensional pointed Hopf
    algebras}, Ann. Math. 171 (2010) 375--417
  [arXiv:\linebreak[0]math$/$0502157 [math.QA]].

\bibitem{[AS-pointed]}N.\;Andruskiewitsch and H.-J.\;Schneider,
  \textit{Pointed Hopf algebras},
  in: \textsl{New directions in Hopf algebras}, MSRI Publications 43,
  pages 1--68.  Cambridge University Press, 2002.

\bibitem{[Andr-remarks]}N.\;Andruskiewitsch, \textit{Some remarks on
    Nichols algebras}, in: \textsl{Hopf algebras}, Bergen, Catoiu and
  Chin (eds.) 25--45. M.\;Dekker (2004).

\bibitem{[Heck-class]}I.\;Heckenberger, \textit{Classification of
    arithmetic root systems}, Adv. Math. 220 (2009) 59--124
  [math.QA$/$\linebreak[0]0605795].

\bibitem{[Heck-Weyl]}I.\;Heckenberger, \textit{The Weyl groupoid of a
    Nichols algebra of diagonal type}, Invent. Math. 164 (2006)
  175--188.

\bibitem{[AHS]}N.\;Andruskiewitsch, I.\;Heckenberger, and
  H.-J.\;Schneider, \textit{The Nichols algebra of a semisimple
    Yetter--Drinfeld module}, Amer. J. Math., 132 (2010) 1493--1547
  [arXiv:\linebreak[0]0803.2430 [math.QA]].

\bibitem{[ARS]}N.\;Andruskiewitsch, D\;Radford, and H.-J.\;Schneider,
  \textit{Complete reducibility theorems for modules over pointed Hopf
    algebras}, J.\ Algebra, 324 (2010) 2932--2970
  [arXiv:\linebreak[0]1001.3977].

\bibitem{[GHV]}M.\;Gra\~na, I.\;Heckenberger, and L.\;Vendramin,
  \textit{Nichols algebras of group type with many quadratic
    relations}, arXiv:\linebreak[0]1004.3723.

\bibitem{[GH-lyndon]}M.\;Gra\~na and I.\;Heckenberger, \textit{On a
    factorization of graded Hopf algebras using Lyndon words}, J.
  Algebra 314 (2007) 324--343.

\bibitem{[AFGV]}N.\;Andruskiewitsch, F.\;Fantino, G.A.\;Garcia, and
  L.\;Vendramin, \textit{On Nichols algebras associated to simple
    racks}, Contemp. Math. 537, Amer. Math. Soc., Providence, RI,
  2011, pp. 31-56.  [arXiv:\linebreak[0]1006.\linebreak[0]5727].

\bibitem{[AAY]}N.\;Andruskiewitsch, I.\;Angiono, and H.\;Yamane,
  \textit{On pointed Hopf superalgebras}, Aparecer\'a en
  Contemp. Math. 544, Amer. Math. Soc., Providence, RI, 2011
  [arXiv:\linebreak[0]1009.5148].

\bibitem{[Ag-0804-standard]}I.\;Angiono, \textit{On Nichols algebras
    with standard braiding}, Algebra \& Number Theory 3 (2009) 35--106
  [arXiv:\linebreak[0]0804.0816].

\bibitem{[Ag-1008-presentation]}I.E.\;Angiono, \textit{A presentation
    by generators and relations of Nichols algebras of diagonal type
    and convex orders on root systems}, arXiv:\linebreak[0]1008.4144.

\bibitem{[Ag-1104-diagonal]}I.\;Angiono, \textit{On Nichols algebras
    of diagonal type}, arXiv:\linebreak[0]1104.0268.

\bibitem{[FHST]} J.\;Fuchs, S.\;Hwang, A.M.\;Semikhatov, and
  I.Yu.\;Tipunin, \textit{Nonsemisimple fusion algebras and the
    Verlinde formula}, Commun. Math. Phys. 247 (2004) 713--742
  [arXiv:\linebreak[0]hep-th/\linebreak[0]0306274].

\bibitem{[AM-3]} D.\;Adamovi\'c and A.\;Milas, \textit{Lattice
    construction of logarithmic modules for certain vertex algebras},
  Selecta Math. New Ser. 15 (2009) 535--561
  [arXiv:\linebreak[0]0902.3417 [math.QA]].

\bibitem{[Kausch]}H.G.\;Kausch, \textit{Extended conformal algebras
    generated by a multiplet of primary fields}, Phys.\ Lett. B~259
  (1991) 448.

\bibitem{[Gaberdiel-K]}M.R.\;Gaberdiel and H.G.\;Kausch,
  \textit{Indecomposable fusion products}, Nucl.\ Phys.\ B477 (1996)
  293--318 [hep-th$/$\linebreak[0]9604026];

\bibitem{[Gaberdiel-K-2]}M.R.\;Gaberdiel and H.G.\;Kausch, \textit{A
    rational logarithmic conformal field theory}, Phys.\ Lett. B~386
  (1996) 131--137 [hep-th$/$\linebreak[0]9606050].

\bibitem{[Gaberdiel-K-3]}M.R.\;Gaberdiel and H.G.\;Kausch, \textit{A
    local logarithmic conformal field theory}, Nucl.\ Phys.\ B538
  (1999) 631--658 [hep-th$/$\linebreak[0]9807091].

\bibitem{[MN]}J.\;Murakami, K.\;Nagatomo, \textit{Logarithmic knot
    invariants arising from restricted quantum groups},
  Internat. J. Math 19 (2008) 1203--1213 [arXiv:0705.3702 [math.GT]].

\bibitem{[FHT]} P.\;Furlan, L.\;Hadjiivanov, and I.\;Todorov,
  \textit{Zero modes' fusion ring and braid group representations for
    the extended chiral WZNW model}, Lett.\ Math.\ Phys.\ 82 (2007)
  117--151 [arXiv:0710.1063].

\bibitem{[Ar]} Y.\;Arike, \textit{Symmetric linear functions of the
    restricted quantum group $\bar{U}_qsl_2(\mathbb{C})$},
  arXiv:\linebreak[0]0706.1113; \textit{Symmetric linear functions on
    the quantum group $g_{p, q}$}, arXiv:\linebreak[0]0904.0331
  [math.QA].

\bibitem{[KoSa]}H.\;Kondo and Y.\;Saito, \textit{Indecomposable
    decomposition of tensor products of modules over the restricted
    quantum universal enveloping algebra associated to
    $\boldsymbol{\mathfrak{sl}_2}$}, arXiv:0901.4221 [math.QA].

\bibitem{[S-yd]}A.M.\;Semikhatov, \textit{Heisenberg double
    $\mathscr{H}(B^*)$ as a braided commutative Yetter--Drinfeld
    module algebra over the Drinfeld double}, Commun. Algebra 39
  (2011) 1883--1906 [arXiv:1001.0733].

\bibitem{[AGL]}A.\;Alekseev, D.\;Gluschenkov, and A.\;Lyakhovskaya,
  \textit{Regular representation of the quantum group $sl_q(2)$
    \textup{(}$q$ is a root of unity\textup{)}}, St.\ Petersburg
  Math. J.\ 6 (1994) 88.

\bibitem{[Su]}R.\;Suter, \textit{Modules over
    $\mathfrak{U}_q(\mathfrak{sl}_2)$}, Commun. Math. Phys. 163 (1994)
  359--393.

\bibitem{[X]}J.\;Xiao, \textit{Finite dimensional representations of
    $U_t(sl(2))$ at roots of unity}, Can. J. Math. 49 (1997) 772--787.


\bibitem{[S-q]} A.M.\;Semikhatov, \textit{Factorizable ribbon quantum
    groups in logarithmic conformal field theories}, Theor.\
  Math. Phys.\ 154 (2008) 433--453 [arXiv:0705.4267 [hep-th]].

\bibitem{[Brug]} A.\;Brugui\`eres, \textit{Double braidings, twists
    and tangle invariants}, J. Pure Appl. Algebra 204 (2006) 170--194.

\bibitem{[Rosso-CR]} M.\;Rosso, \textit{Groupes quantiques et
    alg{\'e}bres de battage quantiques}, C.R.A.S.\ Paris, 320
  (S{\'e}rie~I) (1995) 145--148.

\bibitem{[Sch-borel]} P.\;Schauenburg, \textit{A characterization of
    the Borel-like subalgebras of quantum enveloping algebras},
  Commun. Algebra 24 (1996) 2811--2823.

\bibitem{[A-about]} N.\;Andruskiewitsch, \textit{About finite
    dimensional Hopf algebras},
  Contemp. Math 294 (2002) 1--57.

\bibitem{[G-free]} M.\;Gra\~na, \textit{A freeness theorem for Nichols
    algebras}, J. Alg. 231 (2000) 235--257.

\bibitem{[Rosso-inv]}M.\;Rosso \textit{Quantum groups and quantum
    shuffles}, Invent. math. 133 (1998) 399--416.

\bibitem{[Wor]} S.L.\;Woronowicz, \textit{Differential calculus on
    compact matrix pseudogroups \textup{(}quantum groups\textup{)}},
  Commun. Math. Phys. 122 (1989) 125--170.

\bibitem{[Besp-TMF]}Yu.N.\;Bespalov, \textit{Crossed modules, quantum
    braided groups, and ribbon structures}, Theor.\ Math.\ Phys. 103
  (1995) 621--637.

\bibitem{[Besp-next]}Yu.N.\;Bespalov, \textit{Crossed modules and
    quantum groups in braided categories},
  arXiv:\linebreak[0]q-alg/9510013.

\bibitem{[Besp-Dr-(Bi)]}Yu.\;Bespalov and B.\;Drabant, \textit{Hopf
    (bi-)modules and crossed modules in braided monoidal categories},
  J. Pure and Applied Algebra 123 (1998) 105--129.

\bibitem{[Majid-book]}S.\;Majid, \textit{Foundations of Quantum Group
    Theory}, Cambridge University Press, 1995.

\bibitem{[Sch-H-YD]}P.\;Schauenburg, \textit{Hopf modules and
    Yetter--Drinfel'd modules}, J. Algebra 169 (1994) 874--890.

\bibitem{[2-boson]}A.M.\;Semikhatov, \textit{Virasoro central charges
    for Nichols algebras}, arXiv:\linebreak[0]1109.1767.

\bibitem{[BKLT]} Yu.\;Bespalov, T.\;Kerler, V.\;Lyubashenko, and
  V.\;Turaev, \textit{Integrals for braided Hopf algebras},
  q-alg/\linebreak[0]9709020.
\end{thebibliography}
\end{document}